\documentclass[11pt]{amsart} 
\usepackage[mathscr]{eucal}
\usepackage{amsmath,amsfonts}

\parskip=\smallskipamount

\newtheorem{theorem}{Theorem}[section]

\newtheorem{proposition}[theorem]{Proposition}
\newtheorem{corollary}[theorem]{Corollary}
\newtheorem{lemma}[theorem]{Lemma}

\newtheorem{remark}[theorem]{Remark}

\makeatletter
\@addtoreset{equation}{section}
\makeatother

\newcommand{\CC}{{\mathbb C}}
\newcommand{\NN}{{\mathbb N}}

\newcommand{\FF}{{\mathbb F}}

\newcommand{\cA}{{\mathcal A}}

\newcommand{\cD}{{\mathcal D}}

\newcommand{\cG}{{\mathcal G}}
\newcommand{\cH}{{\mathcal H}}
\newcommand{\cK}{{\mathcal K}}

\newcommand{\cM}{{\mathcal M}}

\newcommand{\cO}{{\mathcal O}}

\newdimen\expt
\def\boxit#1{\setbox0\hbox{$\displaystyle{#1}$}
      \hbox{\lower.4\expt
 \hbox{\lower3\expt\hbox{\lower\dp0
      \hbox{\vbox{\hrule height.4\expt
 \hbox{\vrule width.4\expt\hskip3\expt
      \vbox{\vskip3\expt\box0\vskip2\expt}%
 \hskip3\expt\vrule width.4\expt}\hrule height.4\expt}}}}}}
 
\begin{document}

 


\title [ Similarity and ergodic theory of positive linear maps ] 
{ Similarity and ergodic theory of  positive linear  maps  } 
 \author{Gelu Popescu}
\date{May 1, 2002 (revised version: October  2002)}
\thanks{The author was
partially supported by an NSF grant}

\address{Department of Mathematics, The University of Texas 
at San Antonio \\ San Antonio, TX 78249, USA}
\email{\tt gpopescu@math.utsa.edu}

\begin{abstract} 
In this paper we study   the operator inequality
$\varphi(X)\leq X$ and  the operator equation $\varphi(X)= X$, where $\varphi$ 
is a
$w^*$-continuous  positive
  (resp.\,completely positive)
   linear map   on $B(\cH)$. We show that   their solutions are 
   in one-to-one correspondence
    with  a class of Poisson transforms
   on Cuntz-Toeplitz  $C^*$-algebras, if
    $\varphi$ is  completely positive.  
   Canonical decompositions, ergodic type theorems, and lifting 
   theorems are obtained and used to provide a complete description of
    all solutions, when $\varphi(I)\leq I$.
    
    We show that the above-mentioned  inequality 
    (resp.\,equation) and the structure of its
    solutions 
    have strong implications in connection with representations of 
    Cuntz-Toeplitz  $C^*$-algebras,
    common invariant subspaces for $n$-tuples of operators,
   similarity of positive linear maps, and numerical invariants associated with
  Hilbert modules  
 over  $\CC \FF_n^+$, 
  the complex free semigroup algebra 
   generated by the free semigroup on $n$ generators. 
\end{abstract}

\maketitle

\section{Introduction}\label{INTR}

 Let $\cH$ be a separable Hilbert space and $B(\cH)$ be
  the algebra of all bounded
linear operators on $\cH$. 
  Given 
    a positive linear map $\varphi:B(\cH)\to B(\cH)$,   we define 
         the following sets:
       \begin{enumerate}
\item [(i)]
             $
       C_\leq(\varphi)^+:=\{X\in B(\cH): X\geq 0 \text{ and } 
       \ \varphi(X)\leq X \}
       $ (noncommutative cone);
\item[(ii)]
             $
       C_=(\varphi):=\{X\in B(\cH): \varphi(X)=X \}
       $ (fixed-point operator space).
\end{enumerate}
We will refer to these sets as the  $C(\varphi)$-sets associated with $\varphi$.
  The structure of these  sets plays  a distinguished
    role in  the ergodic theory of  positive maps  
       \cite{DuS}, 
      the classification
    of the endomorphisms of $B(\cH)$ (eg.\,\cite{Pow1},  \cite{Pow2},
    \cite{Arv2}, \cite{Arv3}, \cite{L}, \cite{BJP}, \cite{BrJ1}, \cite{BJKW}),
    and the representation theory of  Cuntz algebras
       (eg.\,\cite{Cu},  \cite{Po-isometric}, \cite{L}, \cite{BJP},
        \cite{BrJ1}, \cite{BJKW}, \cite{DP}, \cite{DKS}).
      In the particular  case when   $T\in B(\cH), \|T\|\leq 1$, 
      and $\varphi_T(X):=TXT^*$
       these sets
      were studied by R.G.~Douglas in \cite{Do} and by Sz.-Nagy 
      and Foia\c s \cite{SzF1}
      in connection with $T$-Toeplitz operators. 
  The operator space $ C_=(\varphi_T)$ was also  studied in 
  \cite{CF1} and \cite{CF2}.
  
  \smallskip
  
  In  this paper we study the structure of the 
  $C(\varphi)$-sets associated with  a $w^*$-continuous  positive
  (resp. completely positive)
   linear map $\varphi$ on $B(\cH)$ and its connections with Poisson transforms
   on Cuntz-Toeplitz  $C^*$-algebras, common invariant subspaces for $n$-tuples of operators,
   similarity of positive linear maps, and numerical invariants associated with
  Hilbert modules  
 over  $\CC \FF_n^+$, 
  the complex free semigroup algebra 
   generated by the free semigroup on $n$ generators.  
  
  It is well-known (see eg.\,\cite{E}) that 
 any $w^*$-continuous completely positive linear map  $\varphi$ 
on $B(\cH)$  is determined by a sequence
 $\{A_i\}_{i=1}^n$ ~$(n\in \NN$ or $n=\infty)$
 of bounded operators on $\cH$, in the sense that
\begin{equation}\label{phiA}
\varphi(X):= \sum_{i=1}^n A_i X A_i^*, \quad X\in B(\cH),
\end{equation}
where, if $n=\infty$, the convergence is in the $w^*$-topology.
 In Section \ref{Poisson}, we show  that 
   the positive solutions of the operator inequality $\varphi(X)\leq X$ are
 intimately related to  a class of  Poisson tranforms on 
 $C^*(S_1,\ldots, S_n)$, the Cuntz-Toeplitz $C^*$-algebra generated by
the left creation operators $S_1,\ldots, S_n$ on the full Fock space.
 More precisely, we  prove  that 
 an operator   $D$ is in $ C_\leq(\varphi)^+$ if and only if  there is 
   a Poisson transform
$$
P_{\varphi, D}: C^*(S_1,\ldots, S_n)\to B(\cH) 
 $$
  with 
   the following properties:
 \begin{enumerate}
 \item[(i)] $P_{\varphi, D}$ is a completely positive linear map;
 \item[(ii)] $\|P_{\varphi, D}\|_{cb}\leq \|D\|$;
 \item[(iii)] $P_{\varphi, D} (I)=D$ and 
 $$ P_{\varphi, D}(S_\alpha S_\beta^*) =A_\alpha D A_\beta^*, 
 \quad \alpha, \beta \in \FF_n^+.
 $$
 \end{enumerate}
When $A_iA_j=A_jA_i$, $i,j=1,\ldots, n$,  the result remains true if we replace 
the  left creation operators $S_1,\ldots, S_n$ by their compressions
$B_1,\ldots, B_n$ to the symmetric Fock space.

Let us  mention that, in the particular case when $\varphi(I)\leq I$ and $D:=I$,
the  Poisson transform  associated with $(\varphi, I)$
was  introduced and studied
in \cite{Po-poisson} in connection with a noncommutative von Neumann inequality
for row contractions \cite{Po-von}. 
Several applications of these   Poisson transforms
 were considered in  \cite{Po-poisson}, \cite{Arv4},
\cite{Po-curvature}, \cite{ArPo1}, \cite{ArPo2}, \cite{Po-tensor}, 
and recently in 
 \cite{BB}, \cite{BBD},  \cite{Po-moment},
and
\cite{Ar}.  We refer to  \cite{Arv1},  \cite{P-book}, and  \cite{Pi-book} for 
results on completely bounded maps and operator spaces.

 In Section \ref{Ergodic}, we present  canonical decompositions (see Theorem
  \ref{decomp}), ergodic type results (see Theorem \ref{limi}),
  and lifting theorems (see Theorem \ref{ineq-pure})  for 
   $w^*$-continuous  positive linear maps   on $B(\cH)$.
   These results together with those from Section \ref{Poisson} are used to provide  a complete description
   of the $C(\varphi)$-sets (see Theorem \ref{main1} and Corollary \ref{cons1}),
   when $\varphi$ is a $w^*$-continuous 
   completely positive linear map with $\varphi(I)\leq I$.
   When we drop the condition  $\varphi(I)\leq I$, we also 
   obtain characterizarions of  the $C(\varphi)$-sets (see Theorem \ref{pure} 
    and Theorem \ref{pure2}).
   An important role in this investigation is played by the noncommutative 
   dilation theory
   for sequences of operators \cite{Fr}, \cite{Bu}, 
   \cite{Po-models}, \cite{Po-isometric}, \cite{Po-charact},
   and \cite{Po-intert} (see \cite{SzF-book} for the classical dilation theory).
  For related results when $\varphi(I)=I$ we mention \cite{BJKW}.

  In Section \ref{Invariant}, we show that there is a strong connection between 
  the positive solutions of the operator inequality 
  $\varphi(X)\leq X$, where $\varphi$ is 
  a $w^*$-continuous completely positive linear map on $B(\cH)$, 
  defined as in \eqref{phiA},
  and the common invariant subspaces for the $n$-tuple of
   operators $\{A_i\}_{i=1}^n$.
   In this direction, we obtain invariant subspace theorems (eg.\,Theorem 
   \ref{fix}) and  Wold type  decomposition theorems for  
   $w^*$-continuous completely  positive linear maps on $B(\cH)$   
   (eg.\,Theorem \ref{wold1}).
  The latter results generalize 
    the classical
    Wold decomposition for isometries, as well as the one obtained
     in \cite{Po-isometric}
    for isometries with orthogonal ranges.
  
  Section \ref{Similarity} is devoted to similarity of positive linear maps on 
  $B(\cH)$.
We say that two linear maps $\varphi, \lambda:B(\cH)\to B(\cH)$ are similar
if there is an invertible operator $R\in B(\cH)$ such that
\begin{equation*} 
\varphi(RXR^*)= R\lambda(X)R^*, \quad \text{ for any } X\in B(\cH).
\end{equation*}
 Notice that this relation is equivalent to
\begin{equation*}
\varphi = \psi_R\circ\lambda \circ \psi_R^{-1},
\end{equation*}
where $\psi_R(X):= RXR^*$, \ $X\in B(\cH)$.
This shows that the  discrete semigroups of completely positive maps
$\{\varphi^k\}_{k=0}^\infty$ and $\{\lambda^k\}_{k=0}^\infty$ are also  similar.
Moreover, $D\in C_\leq (\lambda)^+$ if and only if 
$RDR^*\in C_\leq (\varphi)^+$.
In this section we provide necessary and sufficient conditions 
for a $w^*$-continuous positive linear map $\varphi$ on $B(\cH)$ to be
similar to a positive  linear map $\lambda$ on $B(\cH)$ satisfying one of the following
properties:
\begin{enumerate}
\item[(i)]
$\lambda(I)=I$ (see Theorem \ref{simi});
\item[(ii)] 
$\|\lambda\|<1 $ (see Theorem \ref{simi2});
\item[(iii)] $\lambda$ is a pure completely positive linear map with
$\|\lambda\|\leq 1 $ (see Theorem \ref{sim-pure});
\item[(iv)]
$\lambda$ is a  completely positive linear map with
$\|\lambda\|\leq 1 $ (see Theorem \ref{simi4}).
\end{enumerate}
We show that these similarities are strongly related to the existence of 
invertible positive solutions of the operator inequality
$\varphi(X)\leq X$ or equation  $\varphi(X)= X$.

  In \cite{Arv5}, Arveson 
    introduced  a notion of 
        curvature   and Euler characteristic  for finite rank
       contractive Hilbert modules over $\CC[z_1,\ldots, z_n]$, 
       the complex unital algebra of all polynomials in $n$ commuting variables.
    Noncommutative analogues of these notions were introduced  and studied 
    by the author
   in  \cite{Po-curvature} and, independently, by D.~Kribs  \cite{Kr}.
 The Poisson transforms of Section \ref{Poisson} are used in Section
 \ref{Invariants}  to define certain  numerical  invariants associated 
 with  (not necessarily contractive) Hilbert modules  
 over  
  the  free semigroup algebra  $\CC \FF_n^+$. 
  We extend and refine some of the results from \cite{Po-curvature}.
  Any Hilbert module  $\cH$ over $\CC\FF_n^+$ corresponds to a unique $w^*$-continuous
   completely positive map $\varphi$ on $B(\cH)$ and therefore to a unique noncommutative cone
   $C_\leq (\varphi)^+$.
   A notion of $*$-curvature $\text{ curv}_*(\varphi, D)$ 
   and Euler characteristic $\chi(\varphi, D)$ are associated with 
   each ordered pair $(\varphi, D)$, where $D\in C_\leq (\varphi)^+$.
  In this section, 
  we  obtain  asymptotic 
   formulas and basic properties for both the $*$-curvature 
    and the Euler characteristic associated with $(\varphi, D)$.  
 In the particular case when $\cH$ is a
   contractive Hilbert modules over $\CC\FF_n^+$ and $D:=I$, our two variable invariant
   $$
   F(\varphi, I):=(\|\varphi^*(I)\|, \text{ curv}_*(\varphi, I))
   $$ 
   is a refinement
   of the curvature invariant  from \cite{Po-curvature}  and \cite{Kr}.

\bigskip

\section{ Poisson transforms associated with  completely positive
maps }\label{Poisson}

A Poisson transform on the Cuntz-Toeplitz algebra  
$C^*(S_1,\ldots, S_n)$ is associated with each pair
 $(\varphi, D)$, where $\varphi$ is a 
 $w^*$-continuous completely positive linear map on $B(\cH)$ and $D\in B(\cH)$ 
 is a positive operator such that $\varphi(D)\leq D$.
 The main result of this section (see Theorem \ref{poisson}) shows that 
 the elements of the 
noncommutative cone $C_\leq (\varphi)^+$ are in one-to-one correspondence
 with the elements of a class 
of
Poisson transforms on  Cuntz-Toeplitz algebras.
On the other hand, we prove that 
there is a strong connection  between the 
fixed-point operator space $C_=(\varphi)$ and a class of Poisson transforms
on the Cuntz algebra $\cO_n$.

Let $H_n$ be an $n$-dimensional complex  Hilbert space with orthonormal basis
$e_1,e_2,\dots,e_n$, where $n\in \{1,2,\dots\}$ or $n=\infty$.
  We consider the full Fock space  of $H_n$ defined by
$$F^2(H_n):=\bigoplus_{k\geq 0} H_n^{\otimes k},$$ 
where $H_n^{\otimes 0}:=\CC 1$ and $H_n^{\otimes k}$ is the (Hilbert)
tensor product of $k$ copies of $H_n$.
Define the left creation 
operators $S_i:F^2(H_n)\to F^2(H_n), \  i=1,\dots, n$,  by
$$
 S_i f:=e_i\otimes f , \ f\in F^2(H_n).
$$
 The noncommutative analytic Toeplitz algebra $F_n^\infty$ is the
  WOT-closed algebra generated by the left creation operators
   $S_1,\dots, S_n$  
    and the identity.
 This algebra   and  its norm-closed version (the noncommutative disc
 algebra  $\cA_n$)  were introduced by the author  in \cite{Po-von}
  in connection
   with a noncommutative von Neumann inequality.

Let  $\FF_n^+$ be the free semigroup with $n$ generators $g_1,\dots, g_n$ and
  neutral element $g_0$. 
The length of $\alpha\in\FF_n^+$ is defined by
$|\alpha|:=k$, if $\alpha=g_{i_1}g_{i_2}\cdots g_{i_k}$, and
$|\alpha|:=0$, if $\alpha=g_0$.
We also define
$e_\alpha :=  e_{i_1}\otimes e_{i_2}\otimes \cdots \otimes e_{i_k}$  
 and $e_{g_0}:= 1$.
It is  clear that 
$\{e_\alpha:\alpha\in\FF_n^+\}$ is an orthonormal basis of $F^2(H_n)$.
 If $T_1,\dots,T_n\in B(\cH)$, define 
$T_\alpha :=  T_{i_1}T_{i_2}\cdots T_{i_k}$,
if $\alpha=g_{i_1}g_{i_2}\cdots g_{i_k}$ and 
$T_{g_0}:=I$, the identity on $\cH$.  

Let  $\cH$ be  a separable Hilbert space.  
 Any $w^*$-continuous completely positive linear map  $\varphi$ 
on $B(\cH)$  is determined by a sequence
 $\{A_i\}_{i=1}^n$ ~$(n\in \NN$ or $n=\infty)$
 of bounded operators on $\cH$, in the sense that
$$
\varphi(X):= \sum_{i=1}^n A_i X A_i^*, \quad X\in B(\cH),
$$
where, if $n=\infty$, the convergence is in the $w^*$-topology. 
 Fix such a map $\varphi$ and  
let 
$D\in B(\cH)$ be a positive operator such that $\varphi(D)\leq D$.
Denote $\varphi_r:= r^2 \varphi$, \ $0<r\leq 1$, and define the defect operator
$\Delta_r:=[D-\varphi_r(D)]^{1/2}$.
Notice that, if $0<r<1$, then 
\begin{equation*}\begin{split}
\sum_{k=0}^\infty \varphi_r^k(\Delta_r^2)&=
D-\varphi_r(D)+ \varphi_r(D-\varphi_r(D))+\cdots\\
&= D-\lim_{n\to\infty} r^{2n}\varphi^n(D)= D.
\end{split}
\end{equation*}
If $r=1$, then $\sum_{k=0}^\infty \varphi_r^k(\Delta_r^2)=
D-\varphi^\infty(D)$,
 where  
$\varphi^\infty(D):= \text{\rm SOT}-\lim\limits_{k\to \infty} \varphi^k(D)$
 exists since the  sequence of positive operators
  $\{\varphi^k(D)\}_{k=0}^\infty$
  is decreasing.

We introduce the Poisson kernel associated with 
the ordered pair
$(\varphi, D)$ as the family 
  of operators $K_{\varphi, D, r}: \cH\to F^2(H_n)\otimes \cH$, ~$0<r\leq 1$, 
defined by
\begin{equation}\label{kernel}
K_{\varphi, D, r}h:=\sum_{k=0}^\infty \sum_{\alpha\in \FF_n^+, |\alpha|=k}
e_\alpha\otimes r^{|\alpha|} 
\Delta_r A_\alpha^*h, \quad h\in \cH.
\end{equation}
When $r=1$, we denote $K_{\varphi, D}:=K_{\varphi, D, 1}$.
Notice that, if $0<r<1$, then
\begin{equation}\label{K*K}
K_{\varphi, D, r}^* K_{\varphi, D, r}=\sum_{k=0}^\infty \varphi_r^k(\Delta_r^2)=D.
\end{equation}
When $r=1$, we have
\begin{equation}\label{r=1}
K_{\varphi, D}^* K_{\varphi, D}=D-\varphi^\infty(D).
\end{equation}
 Due to relation  \eqref{kernel},
  for any $ i=1,\ldots, n$,  and $0<r\leq 1$, 
  we have   
 \begin{equation}\label{KASK}
 K_{\varphi, D, r}(rA_i^*)=(S_i^*\otimes I)K_{\varphi, D, r}.
 \end{equation}
 Let $C^*(S_1,\dots,S_n)$ be the
$C^*$-algebra generated by $S_1,\dots,S_n$.
 For  $0<r\leq 1$, define the operator $P_{\varphi, D, r} :C^*(S_1,\ldots,S_n) \to B(\cH)$
 by setting
 \begin{equation}\label{pois}
 P_{\varphi, D, r}(f):= K_{\varphi, D, r}^* (f\otimes I) K_{\varphi, D, r}, 
 \quad f\in C^*(S_1,\ldots, S_n).
 \end{equation}
  Using relation \eqref{KASK} when   $0<r< 1$, we have
\begin{equation}\label{KSK}
K_{\varphi, D, r}^* (S_\alpha S_\beta^*\otimes I) K_{\varphi, D, r}=
r^{|\alpha|+ |\beta|} A_\alpha DA_\beta,\quad \alpha,\beta\in \FF_n^+.
\end{equation}
Hence,  and using relations \eqref{K*K} and  \eqref{pois}, we infer that 
$P_{\varphi, D, r}$   is a completely positive  linear map and 
\begin{equation}\label{vonn}
\|P_{\varphi, D, r}\|_{cb}\leq \|D\|, \quad \text{ for any } 0<r<1.
\end{equation}

\smallskip

Now we can prove the main result of this section, which shows that the elements of the 
noncommutative cone $C_\leq (\varphi)^+$ are in one-to-one correspondence
 with the elements of a class 
of
Poisson transforms on  Cuntz-Toeplitz algebras.
\begin{theorem}\label{poisson}
Let $\varphi$ be a $w^*$-continuous completely positive linear map
 on $B(\cH)$ defined by
$$
\varphi(X):= \sum_{i=1}^n A_i X A_i^*, \quad X\in B(\cH),
$$
and let 
$D\in B(\cH)$ be a positive operator such that $\varphi(D)\leq D$. Then 
the Poisson transform 
$$
P_{\varphi, D}: C^*(S_1,\ldots, S_n)\to B(\cH), \quad 
 P_{\varphi, D}(f):=\lim_{r \to 1} 
 K_{\varphi, D, r}^* (f\otimes I) K_{\varphi, D, r},
 $$
 where the limit exists in the uniform norm,
  has the following properties:
 \begin{enumerate}
 \item[(i)] $P_{\varphi, D}$ is a completely positive linear map;
 \item[(ii)] $\|P_{\varphi, D}\|_{cb}\leq \|D\|$;
 \item[(iii)] $P_{\varphi, D} (I)=D$ and 
 $$ P_{\varphi, D}(S_\alpha S_\beta^*) =A_\alpha D A_\beta^*, 
 \quad \alpha, \beta \in \FF_n^+.
 $$
 \end{enumerate}
\end{theorem}
\begin{proof}

If $q(S_1,\ldots, S_n):= 
\sum_{\alpha,\beta\in \FF_n^+} a_{\alpha \beta}S_\alpha S_\beta^*$
is a polynomial in $C^*(S_1,\ldots, S_n)$ define
\begin{equation*}
q^D(A_1,\ldots, A_n):=  
\sum_{\alpha,\beta\in \FF_n^+} a_{\alpha \beta}A_\alpha DA_\beta^*.
\end{equation*}
The definition is correct since,
according to \eqref{K*K} and \eqref{KSK}, we have
\begin{equation}\label{von2}
\|q^D(A_1,\ldots, A_n)\|\leq \|D\|\|q(S_1,\ldots, S_n)\|.
\end{equation}
Now, if $f\in C^*(S_1,\ldots, S_n)$ and 
 $q_k(S_1,\ldots, S_n)$ is an arbitrary sequence of polynomials in 
 $C^*(S_1,\ldots, S_n)$ convergent to $f$,
 we define the operator
 \begin{equation}\label{fd}
f^D(A_1,\ldots, A_n):= \lim_{k\to\infty}q_k^D(A_1,\ldots, A_n).
\end{equation}
 Taking into account relation \eqref{von2}, it is clear that the
  operator $f^D$ is well-defined and 
 \begin{equation*} 
\|f^D(A_1,\ldots, A_n)\|\leq \|D\|\|f\|.
\end{equation*}
According to relations \eqref{KSK} and \eqref{vonn}, we have
$$
\|q_k^D(rA_1,\ldots, rA_n)\|\leq \|D\|\|q_k(S_1,\ldots, S_n)\|,
$$
for any $0<r\leq 1$.
 Since $P_{\varphi, D, r}$ is a bounded linear operator, we have
 \begin{equation}\label{fdr}\begin{split}
 f^D(rA_1,\ldots, rA_n)&:= \lim_{k\to\infty} q_k^D(rA_1,\ldots, rA_n)\\
 &=\lim_{k\to\infty}P_{\varphi, D, r}(q_k(S_1,\ldots, S_n))=P_{\varphi, D, r}(f),
 \end{split}
 \end{equation}
  for any $0<r<1$.
 Using  relations \eqref{fd}, \eqref{fdr}, the fact that
  $\|f-q_k\|\to 0$ as $k\to \infty$, and 
 $$
 \lim_{r\to 1}q_k^D(rA_1,\ldots, rA_n)= q_k^D(A_1,\ldots, A_n),
 $$
  we can easily prove that 
 $$
 \lim_{r \to 1} P_{\varphi, D, r}(f)=f^D(A_1,\ldots, A_n)
 $$
 in the uniform norm.
 For any $0<r<1$,  $P_{\varphi, D, r}$ is a completely positive linear  map. 
  Hence, and using 
  relations \eqref{KSK},
 \eqref{vonn},  we infer  that 
 $P_{\varphi, D}$
 is a completely positive map with $\|P_{\varphi, D}\|_{cb}\leq \|D\|$.
  The condition (iii) is clearly satisfied  due to relations \eqref{K*K}
  and \eqref{KSK}. 
 The proof is complete.
\end{proof}

Let us remark that if
$\varphi$ is a pure completely positive map, i.e., 
 $\varphi^k(I)\to 0$ strongly,  as $k\to\infty$, 
then the Poisson transform $P_{\varphi, D}$  satisfies the equation
$$
P_{\varphi, D}f= K^*_{\varphi, D}(f\otimes I)K_{\varphi, D},
\quad f\in C^*(S_1,\ldots, S_n).
$$
We should also mention that Theorem \ref{poisson} actually shows 
 that given a $w^*$-continuous completely positive linear map
$\varphi$  on $B(\cH)$, there exists $D\geq 0$ such that $\varphi(D)\leq D$ 
if and only if there is a 
Poisson transform  $P_{\varphi, D}$ with the  properties (i), (ii), and (iii)
from Theorem 
\ref{poisson}. It remains to prove one implication.
 Indeed, if we assume that $P_{\varphi, D}$ satisfies the above-mentioned
conditions, then 
$$
\varphi(D)=\sum_{i=1}^n A_i DA_i^*=P_{\varphi, D} \left(\sum_{i=1}^n S_iS_i^*
\right)\leq 
P_{\varphi, D}(I)= D.
$$

\begin{corollary}\label{Poisson-commut}
If 
$\varphi(X):= \sum_{i=1}^n A_i X A_i^*$ ~$(n\in \NN$ or $n=\infty)$ is a
$w^*$-continuous completely positive
 linear map
 with $A_i A_j=A_j A_i$, $i,j=1,\ldots, n$, then Theorem $\ref{poisson}$
  remains true
 if we replace the left creation operators $\{S_1,\ldots, S_n\}$ by
  their compressions $\{B_1,\ldots, B_n\}$
 to the symmetric Fock space $F_s^2(H_n)$.
\end{corollary}

\begin{proof}
Since $F_s^2(H_n)\subset F^2(H_n)$ is an invariant subspace under 
each  $S_i^*$, \ $i=1,\ldots, n$, we have
$$
P_{F_s^2(H_n)} S_\alpha S_\beta^*|F_s^2(H_n)= B_\alpha B_\beta^*, 
\quad \alpha, \beta\in \FF_n^+.
$$
On the other hand, since the operators $A_i$ are commuting, the Poisson kernel
 $K_{\varphi, D, r}$ takes values 
in $F_s^2(H_n)\otimes \cH$ for any $0<r<1$.
Hence, and using relation \eqref{KSK}, we deduce that
\begin{equation*}\begin{split}
K_{\varphi, D, r}^* (B_\alpha B_\beta^*\otimes I) K_{\varphi, D, r}&=
K_{\varphi, D, r}^*  (S_\alpha S_\beta^*\otimes I) K_{\varphi, D, r}\\
 &=
r^{|\alpha|+ |\beta|} A_\alpha DA_\beta,\quad \alpha,\beta\in \FF_n^+.
\end{split}
\end{equation*}
The rest of the proof is similar to that of Theorem \ref{poisson}.
\end{proof}

We recall \cite{Cu} that if $n\geq 2$, the Cuntz algebra $\cO_n$ is
 the universal $C^*$-algebra
 generated by elements $v_1,\ldots, v_n$ subject to the relations
 $$
 v_i^* v_j=\delta_{ij} I \ \text{ and } \ \sum_{i=1}^n v_iv_i^*=I.
 $$ 
The following result shows that there is a strong connection  between the 
fixed-point operator space $C_=(\varphi)$ and a class of Poisson transforms
on the Cuntz algebra $\cO_n$. The proof is based on 
 noncommutative dilation theory
\cite{Po-isometric}.

\begin{theorem}\label{Poisson-Cuntz}
If 
$\varphi(X):= \sum_{i=1}^n A_i X A_i^*$  ~$(n\geq 2$ or $n=\infty)$ 
is a $w^*$-continuous completely positive
 linear map on $B(\cH)$ and $D\in B(\cH)$ is an  invertible positive
  solution of the equation
  $\varphi(X)=X$, 
 then there is a unique completely positive  linear map 
  $\Phi_{\varphi, D}:\cO_n \to B(\cH)$
 such that   ~$\Phi_{\varphi, D}(I)= D$
  and 
 $$
 \Phi_{\varphi, D}(v_\alpha v_\beta^*) =A_\alpha D A_\beta^*, 
 \quad \alpha, \beta \in \FF_n^+,
 $$
 where $\{v_1,\ldots, v_n\}$ is a system of generators for the Cuntz algebra $\cO_n$.
 Moreover, if $\varphi(I)=I$ and $D$ is  a positive  operator such that  
$\varphi(D)=D$, then the result remains true.
\end{theorem}

\begin{proof}
Assume $D$ is an invertible positive operator with  $\varphi(D)=D$ and set
 $T_i:= D^{-1/2} A_i D^{1/2}$,
\ $i=1,\ldots, n$. Notice that 
$$
\sum_{i=1}^n T_iT_i^*=D^{-1/2}
 \varphi(D) D^{1/2}=I_\cH.
 $$
 According to \cite{Po-isometric}, the minimal isometric dilation 
 of $[T_1,\ldots,T_n]$ is $[V_1,\ldots, V_n]$, where $V_i$ are isometries 
 on a Hilbert space 
 $\cK\supseteq \cH$ such that 
 $$
 \sum_{i=1}^n V_iV_i^*=I_\cK, \ V_i^*|\cH=T_i^*, \text{ and }
 \bigvee\limits_{\alpha\in \FF_n^+}
 V_\alpha \cH=\cK.
 $$
  Therefore, there is a unique unital completely contractive linear map
  $\Psi:C^*(V_1,\ldots, V_n)\to B(\cH)$ such that $\Psi(V_\alpha V_\beta^*)=
  T_\alpha T_\beta^*$, \ $\alpha, \beta\in \FF_n^+$.
  Hence, we infer that  $\Psi_{\varphi, D}:C^*(V_1,\ldots, V_n)\to B(\cH)$,
   given by
  $\Psi_{\varphi, D}(X):= D^{1/2} \Psi(X) D^{1/2}$, is a completely positive linear map
  such that 
 ~$\Psi_{\varphi, D}(I)= D$
  and 
 $$
 \Psi_{\varphi, D}(V_\alpha V_\beta^*)= D^{1/2} T_\alpha T_\beta^* D^{1/2} =A_\alpha D A_\beta^*, 
 \quad \alpha, \beta \in \FF_n^+.
 $$
Therefore, the map $\Phi_{\varphi, D}$ has the required properties and
$\|\Phi_{\varphi, D}\|_{cb}\leq \|D\|$.

Now, let us assume that $\varphi(I)=I$ and $D$ is only a positive operator
 such that $\varphi(D)=D$. Hence,  if $\epsilon>0$, then  $D+\epsilon I$ 
 is positive invertible
 and $\varphi(D+\epsilon I)=D+\epsilon I$.
 Applying the first part of the theorem, we find a completely positive linear map
 $\Psi_{\epsilon}:\cO_n\to B(\cH)$ such that
 \begin{equation}\label{eps}
 \Psi_\epsilon(v_\alpha v_\beta^*)= 
 A_\alpha DA_\beta +\epsilon A_\alpha A_\beta, \quad \alpha,
  \beta \in \FF_n^+,
 \end{equation}
and 
\begin{equation}\label{enorm}
\|\Psi_\epsilon\|_{cb}\leq \|D+\epsilon I\|.
\end{equation}

If $q(v_1,\ldots, v_n):= 
\sum_{\alpha,\beta\in \FF_n^+} a_{\alpha \beta}v_\alpha v_\beta^*$
is a polynomial in $\cO_n$, define
\begin{equation*}
q^D(A_1,\ldots, A_n):=  
\sum_{\alpha,\beta\in \FF_n^+} a_{\alpha \beta}A_\alpha DA_\beta^*.
\end{equation*}
The definition is correct since,
according to \eqref{enorm}, we have
\begin{equation}\label{von-cunt}
\|q^D(A_1,\ldots, A_n)\|\leq \|D\|\|q(V_1,\ldots, V_n)\|.
\end{equation}
Define $\Phi_{\varphi, D}: \cO_n \to B(\cH)$ by  setting
~$
\Phi_{\varphi, D}f:=f^D(A_1,\ldots, A_n)$, \  $f\in \cO_n$,
where 
\begin{equation}\label{fd-cu}
f^D(A_1,\ldots, A_n):= \lim_{k\to\infty} q_k^D(A_1,\ldots, A_n)
\end{equation}
and 
 $q_k(v_1,\ldots, v_n)$ is an arbitrary sequence of polynomials in 
 $\cO_n$ convergent to $f$.
  Using relation 
\eqref{eps} and standard approximation arguments, we deduce
$$
\Psi_\epsilon(f)= \Phi_{\varphi, D}(f)+\epsilon P_{\varphi, I}(f), 
\quad f\in \cO_n,
$$
where $P_{\varphi, I}$ is the Poisson transform associated 
with $\varphi$ and $I$. Taking $\epsilon\to 0$, we infer that
$$
\Phi_{\varphi, D}(f)=\lim_{\epsilon\to 0} \Psi_\epsilon(f), \quad f\in \cO_n,
$$ in the uniform norm. 
Since  for each $\epsilon>0$,  $\Psi_\epsilon$  is a completely positive
linear map satisfying \eqref{enorm}, we deduce that $\Phi_{\varphi, D}$
 is completely positive 
and $\|\Phi_{\varphi, D}\|_{cb}\leq \|D\|$. The proof is complete.
\end{proof}

\bigskip

\section{ Ergodic theory of completely positive  maps on $B(\cH)$ }
\label{Ergodic}

In  this section  we present  canonical decompositions (see Theorem
  \ref{decomp}), ergodic type results (see Theorem \ref{limi}),
  and lifting theorems (see Theorem \ref{ineq-pure})  for 
   $w^*$-continuous  positive linear maps   on $B(\cH)$.
   These results and the Poisson kernels of Section \ref{Poisson}
     are used  to prove the main result of this section
   (see Theorem \ref{main1}), which 
    provides  a complete description
   of the $C(\varphi)$-sets,
   when $\varphi$ is a $w^*$-continuous 
   completely positive linear map with $\varphi(I)\leq I$.
   When we drop the condition  $\varphi(I)\leq I$, we also 
   obtain characterizations of  the $C(\varphi)$-sets (see Theorem \ref{pure} 
    and Theorem \ref{pure2}).

Let $\varphi$ be a $w^*$-continuous positive linear map on $B(\cH)$. An operator 
$C\in B(\cH)$ is called pure solution of the inequality  $\varphi(X)\leq X$ if
$$\text{  
SOT}-\lim_{k\to \infty} \varphi^k (C)=0.
$$
 Notice that a pure solution is always
a positive operator. In what follows we present a {\it canonical decomposition}
 for the selfadjoint solutions of the operator inequality $\varphi(X)\leq X$.

\begin{theorem}\label{decomp}
Let $\varphi$ be a $w^*$-continuous positive linear map on $B(\cH)$ and
  let $A\in B(\cH)$ be a selfadjoint solution of the inequality
    $\varphi(X)\leq X$.
  Then there exist operators $B,C\in B(\cH)$ with the properties:
  \begin{enumerate}
  \item[(i)] $B=B^*$  is a solution of the equation $\varphi(X)=X$;
  \item[(ii)] $C\geq 0$ is a pure solution  
  of the inequality $\varphi(X)\leq X$;
  \item[(iii)] $A=B+C$.
  \end{enumerate}
   Moreover, this decomposition is unique.
\end{theorem}

\begin{proof}
The sequence  of selfadjoint operators $\{\varphi^k(A)\}_{k=0}^\infty$
 is bounded and decreasing. Thus it converges strongly to to a selfadjoint
  operator $B:= \text{\rm SOT}-\lim\limits_{k\to \infty}\varphi^k(A)$. 
  Since $\varphi$ is $w^*$-continuous,
  we  have $\varphi(B)=B$. Setting $C:=A-B$, we clearly have $C\geq 0$ and 
  $$
  \varphi(C)=\varphi(A)- B\leq A-B=C.
  $$
  Since $\varphi^k(C)\to 0$ strongly, as $k\to \infty$,  $C$ is a pure solution 
   of  the inequality $\varphi(X)\leq X$.
 
 Now suppose $A=B_1+C_1$  with $\varphi(B_1)=B_1$ and $C_1$ 
 is a pure solution  of the inequality  $\varphi(X)\leq X$.
 Then 
 $$
 B-B_1= \varphi^k(B-B_1)=\varphi^k(C_1-C)\to 0, \text{ as } k\to \infty.
 $$
 Therefore,  $B=B_1$ and $C=C_1$. The proof is complete.
\end{proof}
Let us remark that a  result similar to Theorem \ref{decomp} holds
 if  $A$ is a selfadjoint solution
 of the inequality $\varphi(X)\geq X$.

Now, we can prove the following  ergodic type result.
\begin{theorem}\label{limi}
Let $\varphi$ be a w$^*$-continuous  positive linear map on $B(\cH)$
and let $A\in B(\cH)$ be a selfadjoint solution
 of the inequality $\varphi(X)\leq X$.
Then
$$
\text{\rm SOT}-\lim_{k\to\infty} {\frac{\varphi^0(A)+ \varphi^1(A)+\cdots +
 \varphi^{k-1}(A)}{k}}=B,
 $$
 where $A=B+C$ is the canonical decomposition of $A$ with respect to $\varphi$, and 
 $\varphi(B)=B$.
\end{theorem}
\begin{proof}
Since $C$ is a pure solution of the inequality $\varphi(X)\leq X$, we have 
$\text{\rm SOT}-\lim\limits_{k\to\infty}
 \varphi^k(C)=0$. Taking into account that $0\leq\varphi^k(C)\leq C$,
  $k=0,1,\ldots $,
    a standard argument shows that 
 $$
\text{\rm SOT}-\lim_{k\to\infty} {\frac{\varphi^0(C)+ \varphi^1(C)+\cdots +
 \varphi^{k-1}(C)}{k}}=0.
 $$
 On the other hand,
 since $A=B+C$ and $\varphi(B)=B$, we infer that
 $$
 {\frac{\varphi^0(A)+ \varphi^1(A)+\cdots +
 \varphi^{k-1}(A)}{k}}=B+{\frac{\varphi^0(C)+ \varphi^1(C)+\cdots +
 \varphi^{k-1}(C)}{k}}.
 $$
 Hence, the result follows.
\end{proof}

We recall from \cite{Po-isometric} the following  Wold type  decomposition
for isometries with orthogonal ranges. Let $V_i\in B(\cK)$, $i=1,\ldots, n$, be isometries
with $V_i^* V_j=0$ if $i\neq j$. Then there are subspaces $\cK_c, \cK_s\subseteq \cK$ reducing for
each $V_1,\ldots, V_n$, such that 
\begin{enumerate}
\item[(i)] $\cK= \cK_c\oplus \cK_s;$
\item[(ii)] $\left(\sum_{i=1}^n V_i V_i^*\right)| \cK_c= I_{\cK_c};$
\item[(iii)]
$\{V_i|{\cK_s}\}_{i=1}^n$ is unitarily equivalent to 
$\{S_i\otimes I_\cM\}_{i=1}^n$ for some Hilbert space $\cM$.
\end{enumerate}
Moreover, the decomposition is unique, up to a unitary equivalence.
Since the isometries $\{V_i|{\cK_c}\}_{i=1}^n$ generate a 
representation of the Cuntz algebra $\cO_n$, we call $\cK_c$ the Cuntz part
in the Wold decomposition  $\cK= \cK_c\oplus \cK_s$.

\begin{corollary}\label{limcon}
Let $\varphi_T$ be a w$^*$-continuous  completely positive linear map on $B(\cH)$
such that $\varphi_T(I)\leq I$ and $\varphi_T(X)= \sum_{i=1}^n T_i XT_i^*$
 $(n\in \NN$ or $n=\infty)$.
 Then
$$
\text{\rm SOT}-\lim_{k\to\infty} {\frac{\varphi_T^0(I)+ \varphi_T^1(I)+\cdots 
+
 \varphi_T^{k-1}(I)}{k}}=P_\cH P_{\cK_c}|\cH,
 $$
 where $\cK_c$ is the Cuntz part in 
 the Wold decomposition  $\cK=\cK_c\oplus \cK_s$ of the minimal 
 isometric dilation $[V_1,\ldots, V_n]$ of $[T_1,\ldots, T_n]$ on the 
 Hilbert space 
 $\cK\supseteq \cH$.  
\end{corollary}
\begin{proof}
Since $[V_1,\ldots, V_n]$ is the minimal isometric dilation 
of $[T_1, \ldots, T_n]$ on 
$\cK\supseteq \cH$, we have $V_i^*|\cH= T_i^*$, $i=1,\ldots, n$. Define 
the completely positive map
$\varphi_V(Y):= \sum_{i=1}^n V_i Y V_i^*$, \ $Y\in B(\cK)$.
Since $\varphi_T^j(I_\cH)=P_\cH \varphi_V^j (I_\cK)|\cH$, \ $j=1,2,\ldots$, 
it  remains to prove that 
\begin{equation}
\label{erg}
\text{\rm SOT}-\lim\limits_{k\to\infty} {\frac{\varphi_V^0(I)+ 
\varphi_V^1(I)+\cdots 
+
 \varphi_V^{k-1}(I)}{k}}= P_{\cK_c}.
\end{equation}
According to Theorem \ref{limi} and Theorem \ref{decomp}, the limit 
in \eqref{erg} is equal to
$\text{\rm SOT}-\lim\limits_{k\to\infty} \varphi_V^k(I_\cK)$.
 The noncommutative Wold decomposition \cite{Po-isometric}, shows that 
 $$
 \text{\rm SOT}-\lim_{k\to\infty} \varphi_V^k(I_\cK)=P_{\cK_c}.
 $$
 This completes the proof.
\end{proof}

We recall \cite{Po-isometric} that an $n$-tuple of operators 
$[T_1,\ldots, T_n]$, $T_i\in B(\cH)$, is a $C_0$-row contraction if 
$T_1T_1^*+\cdots+T_nT_n^*\leq I_\cH$ and 
$$
\lim_{k\to\infty}  \sum_{\alpha\in  \FF_n^+, |\alpha|=k} \|T_\alpha^* h\|^2=0
\quad \text{ for any } h\in \cH.
$$
In
 what follows,  we  obtain a characterization of the solutions of
 the  inequality    $\varphi(X)\leq X$ (resp.\, equation $\varphi(X)= X$),
 where  
 $\varphi$ is a $w^*$-continuous completely positive linear map on $B(\cH)$.

\begin{theorem}\label{pure}

Let $\varphi$ be a w$^*$-continuous  completely positive linear map on $B(\cH)$
given by
$$
\varphi(X):= \sum_{i=1}^n A_i X A_i^*,\quad X\in B(\cH).
$$
 A positive  operator  $C\in B(\cH)$  is a solution of the inequality 
 $\varphi(X)\leq X$ (resp.\, equation  $\varphi(X)= X$) if and only
  if there exist operators
 $B_i\in B(\cH)$, \ $i=1,\ldots, n$, such that 
 $\sum_{i=1}^n B_iB_i^*\leq 1$ (resp.\,$\sum_{i=1}^n B_iB_i^*= 1$) and 
 \begin{equation}\label{c12}
 A_i C^{1/2}= C^{1/2} B_i, \quad i=1,\ldots, n.
 \end{equation}
 Moreover, $C$ is a pure solution of  $\varphi(X)\leq X$ if and only if 
 there exists a  $C_0$-row contraction $[B_1,\ldots, B_n]$ satisfying
 relation
 \eqref{c12}.
\end{theorem}
\begin{proof}
   Assume that $C\in B(\cH)$  is a solution of 
 the inequality 
 $\varphi(X)\leq X$ (resp.\,equation $\varphi(X)= X$). Define  the operator
 $G_i: \overline{\text{\rm range}~C^{1/2}}\to 
 \overline{\text{\rm range}~C^{1/2}}$  by setting
 $$
 G_i^* C^{1/2}:= C^{1/2} A_i^*, \quad i=1,\ldots, n.
 $$
 The definition is correct since
 \begin{equation}\label{corect}
 \sum_{i=1}^n \|G_i^*C^{1/2}h\|^2 = \sum_{i=1}^n \|C^{1/2} A_i^*h\|^2
 =\langle \varphi(C) h, h\rangle \leq \|C^{1/2}h\|^2.
 \end{equation}
 If $\varphi(C)=C$, then we have equality in \eqref{corect}.
 Let $Q_i$, \ $i=1,\ldots, n$,  
  be bounded operators on 
 $\cM:=(\overline{\text{\rm range}~C^{1/2}})^\perp$ such
  that $\sum_{i=1}^n Q_i Q_i^*=I$.
  Define $B_i:= G_i\oplus Q_i$, $i=1,\ldots, n$,
  with respect to the decomposition
 $\cH= \cM^\perp\oplus \cM$,
   and notice that  
 $\sum_{i=1}^n B_i B_i^*\leq I$ if $\varphi(C)\leq C$, and 
  $\sum_{i=1}^n B_i B_i^*= I$ if $\varphi(C)= C$.
  
  Conversely, assume that $B_i\in B(\cH)$ satisfies
  $
 A_i C^{1/2}= C^{1/2} B_i$,  for any $i=1,\ldots, n$.
 Then  we have
 $$
 \varphi(C)= C^{1/2}\left(\sum_{i=1}^n B_i B_i^*\right) C^{1/2} \leq C
 $$
 if  $\sum\limits_{i=1}^n B_i B_i^*\leq I$, and $\varphi(C)= C$ if
 $\sum\limits_{i=1}^n B_i B_i^*= I$. 
 
 To  prove the second part of the theorem, assume that $C$ is a pure
  solution of $\varphi(X)\leq X$. Following the first part of the proof,
   define $B_i:= G_i\oplus 0$, \ $i=1,\ldots, n$, with respect to the decomposition 
   $\cH=\cM^\perp \oplus \cM$.
   Since 
   $$\sum_{|\alpha|=k}\|G_\alpha^* C^{1/2}h\|^2= \langle \varphi^k(C) h, h\rangle,
   \quad h\in \cH,
   $$
   it is clear that
   $\sum\limits_{|\alpha|=k} B_\alpha B_\alpha^*\to 0$ strongly,  as  $k\to 0$.
   Therefore, $[B_1,\ldots, B_n]$ is a $C_0$-row contraction.
    For the converse,
   it is enough to observe that
   $$
   \varphi^k(C)=C^{1/2} \left( \sum_{|\alpha|=k} B_\alpha B_\alpha^*\right)
    C^{1/2}.
   $$
 This completes the proof. 
\end{proof}

Consider now the case when $\varphi$ is a $w^*$-continuous completely positive 
linear map on $B(\cH)$ with $\varphi(I)\leq I$. 
Let
$$
\varphi_T(X):= \sum_{i=1}^n T_i X T_i^*,\quad X\in B(\cH),
$$
where $n\in \NN$ or $n=\infty$,
 and 
let 
$$
\varphi_V(Y):= \sum_{i=1}^n V_i Y V_i^*,\quad Y\in B(\cK),
$$ 
where $[V_1,\ldots, V_n]$ is the minimal isometric dilation 
of the row contraction $[T_1,\ldots, T_n]$ on the Hilbert
 space $\cK\supseteq \cH$ (see \cite{Po-isometric}), i.e.,
  $V_i^*|_\cH= T_i^*$, \ $i=1,\ldots, n$,  and
 $\cK=\bigvee_{\alpha\in \FF_n^+} V_\alpha \cH$. We call $\varphi_V$ 
 the minimal dilation of 
 $\varphi_T$.  Notice that $\varphi_V$ 
 is a normal $*$-endomorphism of $B(\cH)$ such that
 $$
 \langle \varphi_T^k(X)h,h'\rangle=\langle \varphi_V^k(X)h,h'\rangle,
 $$
 for any $h,h'\in \cH$, \ $X\in B(\cH)\subseteq B(\cK)$, and $k\in \NN$.
 Here we identify $X\in B(\cH)$ with $P_\cH XP_\cH\in B(\cK)$, where
 $P_\cH$ is the orthogonal projection of $\cK$ onto $\cH$.

The noncommutative dilation theory together with Theorem \ref{pure}
can be used to obtain the following lifting theorem for the
 solutions of the  operator inequality $\varphi_T(X)\leq X$ (resp.\,equation
  $\varphi_T(X)= X$).

\begin{theorem}\label{ineq-pure}
Let $\varphi_T$ be a w$^*$-continuous completely positive 
linear map with $\varphi_T(I)\leq I$ and let $\varphi_V$ be its
 minimal isometric dilation. 
 \begin{enumerate}
\item[(i)] A positive operator $C\in B(\cH)$ 
is a solution of the inequality $\varphi_T(X)\leq X$ (resp.\,equation
  $\varphi_T(X)= X$), if and only if 
  $
 C:=P_\cH D|\cH,
 $
 where $D$ is a positive solution of 
 the inequality $\varphi_V(Y)\leq Y$ (resp.\,equation
  $\varphi_V(Y)= Y$),  such that $\|C\|=\|D\|$.
  \item[(ii)] An operator $C\in B(\cH)$ 
is a pure solution of the inequality $\varphi_T(X)\leq X$ if and only if
 $
 C:=P_\cH D|\cH,
 $
where $D$ is a pure solution of the inequality $\varphi_V(Y)\leq Y$,  such that $\|C\|=\|D\|$.
\end{enumerate}
\end{theorem}
\begin{proof}
Assume that $C\in B(\cH)$ 
is a  solution of the inequality $\varphi_T(X)\leq X$.
 Taking into account Theorem \ref{pure}, we  
find $B_i\in B(\cH)$ satisfying
$$
T_i C^{1/2}= C^{1/2} B_i,
\ i=1,\ldots, n,
$$ 
  where $[B_1,\ldots, B_n]$ is a row contraction which has the property
  $\sum_{i=1}^n B_i B_i^*=1$ if $\varphi(C)=I$, and  
  $\sum\limits_{|\alpha|=k} B_\alpha B_\alpha^*\to 0$ strongly,  as $ k\to 0$,
   if $\varphi^k(C)\to 0$.
  Let $[V_1,\ldots, V_n]$ be the minimal isometric dilation of 
   $[T_1,\ldots, T_n]$ on a Hilbert space $\cK_1\supseteq \cH$, and
   let  $[W_1,\ldots, W_n]$ be the minimal isometric dilation 
   of $[B_1,\ldots, B_n]$  on a Hilbert space $\cK_2\supseteq \cH$.
 According to the noncommutative commutant lifting theorem
 \cite{Po-isometric} (see also  \cite{Po-intert}), there exists an operator
  $\tilde{C}:\cK_1\to \cK_2$
 such that $C^{1/2}= \tilde{C}|\cH$, \ $\|\tilde{C}\|=\|C^{1/2}\|$, and 
 $$
 \tilde{C}V_i^*=W_i \tilde{C}, \quad  i=1,\ldots, n.
 $$
 Notice that
 $$
 \varphi_V(\tilde{C}^* \tilde{C})=
  \tilde{C}^* \left(\sum_{i=1}^n W_i W_i^* \right)\tilde{C}\leq \tilde{C}^* \tilde{C},
 $$
 if $\varphi_T(C)\leq C$ and $\varphi_V(\tilde{C}^* 
 \tilde{C})=\tilde{C}^* \tilde{C}$
 if $\varphi_T(C)= C$.
 Setting $D:= \tilde{C}^* \tilde{C}$, we have   \ $\|D\|=\|C\|$,  
  \ $ C= P_\cH D|\cH$, and the first part of the theorem is proved.
  
 For the second part, if $\varphi^k_T(C)\to 0$ strongly, as $k\to \infty$, then,
 according to Theorem \ref{pure},  the row-contraction  $[B_1,\ldots, B_n]$
  is  of class 
  $C_0$  and  its minimal isometric dilation 
  $[W_1,\ldots, W_n]$ is a $C_0$-row isometry.
  Therefore,  
 $$
 \varphi_V^k(D)=
  \tilde{C}^* \left(\sum_{|\alpha|=k} W_\alpha W_\alpha^*\right) 
  \tilde{C}\to 0\ 
   \text{ strongly,  as } k\to \infty.
 $$
 Conversely, if $D$ is a   solution of the inequality $\varphi_V(Y)\leq Y$, 
 then 
 \begin{equation*}\begin{split}
 \varphi(P_\cH D|\cH)&= \sum_{i=1}^n T_i (P_\cH D|\cH)T_i^*=
 P_\cH (\sum _{i=1}^n V_i D V_i^*)|\cH\\
 &= P_\cH \varphi_V(D)|\cH\leq P_\cH D|\cH.
 \end{split}
 \end{equation*}
 Notice that we have equality if $\varphi_V(D)= D$.
 On the other hand,  since  
 $$
 \varphi^k(P_\cH D|\cH)= P_\cH \varphi^k_V(D)|\cH\to 0 \ \text{ strongly, 
 as } k\to \infty,
 $$
  it is clear that $C:=P_\cH D|\cH$ is a pure solution of
  the inequality
 $\varphi_T(X)\leq X$ when $D$ is a pure solution of
  the inequality $\varphi_V(Y)\leq Y$.
  The proof is complete.
\end{proof}

\begin{corollary}\label{brat} \text{\rm (\cite{BJKW})}
If $\varphi_T(I)=I$, then  a positive operator $C\in B(\cH)$ 
is a solution of
  the equation  
$\varphi_T(X)=X$ if and only if there exists 
$D\in \{V_i, V_i^*\}'$ such that
$$
C=P_\cH D|\cH,\ \|D\|=\|C\|.
$$
Moreover, the result remains true if $C$ is a selfadjoint operator.
\end{corollary}
\begin{proof}
Notice that if $\varphi_V(Y)=Y$, then  $Y\in \{V_i, V_i^*\}'$.
 If 
 $\varphi_T(I)=I$, then the converse is  also true.
  In this case we have    $\sum_{i=1}^n V_i V_i^*= I$. 
  Applying Theorem \ref{ineq-pure}, the result follows.
\end{proof}

The Poisson kernels of Section \ref{Poisson} can be used to better understand 
the structure of the pure solutions of the inequality
$\varphi(X)\leq X$, where 
 $\varphi$ is a $w^*$-continuous  completely positive linear map on $B(\cH)$.
\begin{theorem}\label{pure2}
Let $\varphi$ be a $w^*$-continuous  completely positive linear map on $B(\cH)$
given by
$$
\varphi(X):= \sum_{i=1}^n A_i X A_i^*,\quad X\in B(\cH).
$$
 A positive operator $C\in B(\cH)$ is a pure solution of the inequality 
  $\varphi(X)\leq X$
 if and only if there is a Hilbert space $\cD$ and an operator  
 $K:\cH\to F^2(H_n)\otimes \cD$ such that 
 \begin{equation}\label{int}
 C=K^*K\quad \text{ and } \quad KA_i^*= (S_i^*\otimes I_\cD)K, \ i=1,\ldots, n.
 \end{equation}
\end{theorem}

\begin{proof}
Assume $C$ is a positive solution of $\varphi(X)\leq X$. 
Let $K_{\varphi, C}:\cH\to F^2(H_n)\otimes \cD$ be the Poisson kernel associated
 with $\varphi$ and $C$, i.e.,
$$
K_{\varphi, C}h:= \sum_{\alpha\in \FF_n^+} e_\alpha \otimes \Delta A_\alpha^*h, 
\quad h\in \cH,
$$
where $\Delta:= (C-\varphi(C))^{1/2}$ and 
$\cD:=\overline{\text{\rm range}~\Delta}$.
According to the results of Section \ref{Poisson},  the relation 
\eqref{int} holds when $K:=K_{\varphi, C}$ .

Conversely, if  relation \eqref{int} is satisfied, then 
$$
\varphi(K^* K)= \sum_{i=1}^n A_iK^* KA_i^*= K^* 
\left(\sum_{i=1}^n S_iS_i^*\otimes I\right)K\leq K^*K.
$$
Since $\varphi^k(K^* K)= K^*\left(\sum\limits_{|\alpha|=k}^n S_\alpha
 S_\alpha^*\otimes I\right)K$
is strongly convergent to zero as $k\to \infty$, the result follows.
\end{proof}

 The main result of this section is the following   theorem
  which characterizes the positive solutions of the operator inequality 
  $\varphi_T(X)\leq X$, when $\varphi_T(I)\leq I$.
  
 \begin{theorem}\label{main1}
Let $\varphi_T$ be a w$^*$-continuous completely positive 
linear map with $\varphi_T(I)\leq I$ and let $\varphi_V$ be its
 minimal isometric dilation. A positive operator $A\in B(\cH)$ is a   
 solution of the operator inequality 
  $\varphi_T(X)\leq X$
 if and only if there exist 
 \begin{enumerate}
 \item[(i)] an operator $B\in B(\cK)$ 
 in the commutant 
 of $\{V_i, V_i^*\}_{i=1}^\infty$, 
 \item[(ii)]  a Hilbert space $\cD$, and an operator 
 $K:\cK \to F^2(H_n)\otimes \cD$  with 
 $$V_i K^*= K^*(S_i\otimes I_\cD), 
 \quad i=1,2,\ldots,
 $$ 
 \end{enumerate}
  such that
 \begin{equation}\label{mainsol}
 A=P_\cH A_1|_\cH 
 + P_\cH A_2 |_\cH, \quad \|A\|=\|A_1+ A_2\|,
 \end{equation}
 where 
 $$A_1:=[\text{\rm SOT}-\lim_{k\to \infty} \varphi_V^k(B)], \quad 
 A_2:=  K^*K,
 $$
  and 
  $P_\cH$ is the orthogonal projection on $\cH$. 
 Moreover, the canonical decomposition of $A$ with respect to $\varphi_T$ 
 coincides with   the decomposition  from \eqref{mainsol}.
 \end{theorem}
 \begin{proof}
 According to Theorem \ref{ineq-pure},
  a positive operator $A\in B(\cH)$ is a solution of the operator inequality
 $\varphi_T(X)\leq X$ if and only if there exists a positive operator
 $D\in B(\cK)$ such that  
 $$
 \varphi_V(D)\leq D, \ A=P_\cH D|\cH,  \text{ and } \|A\|=\|D\|.
 $$
 Let 
   $D=B+C$ be  the canonical decomposition of $D$  with respect to 
   $\varphi_V$,   i.e.,
   
  \begin{enumerate}
 \item[(a)] $B= B^*$ is a solution of the equation $\varphi_V(Y)=Y$;
 \item[(b)] $C\in B(\cK)$ is a pure solution  of the inequality 
 $\varphi_V(Y)\leq Y$.
 \end{enumerate}
 %
 According to Theorem \ref{decomp}, 
 we have $B=\text{\rm SOT}-\lim\limits_{k\to\infty} \varphi^k_V(D)$.
 Since $\varphi_V(B)= B$,  it is clear that 
 $B$ in  the commutant of $\{V_i, V_i^*\}_{i=1}^\infty$, and $B=\varphi^k_V(B)$, for any 
 $k\in \NN$.
 On the other hand,
 applying  Theorem \ref{pure2} to the operator $C$, we infer that there is a 
   Poisson kernel $K:\cK\to F^2(H_n)\otimes \cD$ associated
 with $\varphi_V$ and $C$, such that
 $$
 KV_i^*=(S_i^*\otimes I) K,\quad  i=1,\ldots, n, 
 $$
 and $ C=K^*K$.
  Summing up and setting $A_2:= C$, we obtain relation \eqref{mainsol}.

 Conversely, if  $B\in B(\cK)$ 
 in the commutant 
 of $\{V_i, V_i^*\}_{i=1}^\infty$, then $\varphi_V(B)\leq B$. Hence 
 $A_1:=[\text{\rm SOT}-\lim\limits_{k\to \infty} \varphi_V^k(B)]$ exists and
 $\varphi_V(A_1)=A_1$. On the other hand, if $K$ satisfies  (ii), then,
  according to 
 Theorem \ref{pure2}, the operator $KK^*$ is a pure solution of the inequality
  $\varphi_V(Y)\leq Y$. Therefore, we have $\varphi_V(A_1+ KK^*)\leq A_1+ KK^*$.
  As in Theorem \ref{ineq-pure}, we infer that
  $$
  A=P_\cH A_1|_\cH 
 + P_\cH KK^* |_\cH
 $$
 is a positive solution of the inequality 
 $\varphi_T(X)\leq X$. The proof is complete.
 \end{proof}

 \begin{corollary}\label{cons1} 
   If $\varphi_T(I)\leq I$ and  $A\in B(\cH)$ is  a positive operator,
   then $\varphi_T(A)=A$ if and only if
 there exists $B\in B(\cK)$, $B\geq 0$,  
 in the commutant 
 of $\{V_i, V_i^*\}_{i=1}^\infty$  such
 $$
  A=P_\cH [\text{\rm SOT}-\lim_{k\to \infty} \varphi_V^k(B)]|_\cH
  =P_\cH B P_{\cK_u} |_\cH,
  $$
  where $\cK_u$ is the Cuntz part in the  
   Wold decomposition  $\cK=\cK_u\oplus \cK_s$ of 
   the minimal isometric dilation $[V_1,\ldots, V_n]$ of $[T_1\ldots, T_n]$ 
   on the Hilbert space $\cK\supseteq \cH$.
 \end{corollary}
  
  Let $\cK$ and $ \cK'$ be Hilbert spaces. We recall from \cite{Po-charact}
  that 
 a bounded linear
  operator 
$M\in B(F^2(H_n)\otimes \cK, F^2(H_n)\otimes \cK')$ is  {\em multi-analytic}
if $M(S_i\otimes I_\cK)= (S_i\otimes I_{\cK'}) M$, \  $i=1,\dots, n$.
The set of multi-analytic operators coincides with  the operator space
$R_n^\infty\bar\otimes B(\cK, \cK')$, where $R_n^\infty$ is the commutant of 
the noncommutative analytic Toeplitz algebra $F_n^\infty$.
More about multi-analytic operators on Fock spaces can be found in
\cite{Po-multi}, 
\cite{Po-analytic}, and \cite{Po-tensor}.
 \begin{corollary}\label{Adv}\text{\rm(\cite{Po-curvature})}
   If $\varphi_T(I)\leq I$ and $\varphi_T $ is pure, i.e., 
   $\varphi_T^k(I)\to 0$ strongly,  as $k\to\infty$,
    then any positive  solution $A$ of  the inequality  $\varphi_T(X)\leq X$  
    is pure and 
    $A=P_\cH \Psi \Psi^*|_\cH $, where $\Psi$ is a multi-analytic operator.
 \end{corollary}
 \begin{proof}
 Since $\varphi_T^k(A)\leq \|A\|\,\varphi_T^k(I)$, $k=1,2, \ldots$,
  and $\varphi_T$ is pure, we infer that 
  any positive  solution  of  the inequality  $\varphi_T(X)\leq X$  
    is pure.
    On the other hand, since $\varphi_T$ is pure, $[T_1,\ldots, T_n]$
    is a $C_0$-row contraction. According to \cite{Po-isometric},
    its minimal isometric dilation ca be identified
    with $[  S_1\otimes I_\cM, \ldots,  S_n\otimes I_\cM]$  for some 
    Hilbert space 
    $\cM$.
    Applying Theorem \ref{main1}, we have 
    $$
    (S_i\otimes I_\cM) K^*= K^* (S_i\otimes I_\cD), \quad i=1,\ldots, n,
    $$
    i.e., $K^*$ is a multi-analytic operator, and the result follows.
 \end{proof}

 \begin{proposition}
If $\varphi_T(I)=I$ and  $A\in B(\cH)$ is a selfadjoint operator, then  
$\varphi_T(A)\leq A$
if and only if 
there exist 
   an operator $B\in B(\cK)$ 
 in the commutant 
 of $\{V_i, V_i^*\}_{i=1}^\infty$, 
  a Hilbert space $\cD$, and an operator 
 $K:\cK \to F^2(H_n)\otimes \cD$  with ~$V_i K^*= K^*(S_i\otimes I_\cD)$, 
 \ $i=1,\ldots, n$, 
  such that
 \begin{equation}\label{mainsol2}
 A=P_\cH B|_\cH 
 + P_\cH KK^* |_\cH, 
 \end{equation}
 where 
  $P_\cH$ is the orthogonal projection on $\cH$. 
 Moreover, the canonical decomposition of $A$ with respect to $\varphi_T$ 
 coincides with   the decomposition    \eqref{mainsol2}.
 \end{proposition}
 \begin{proof}
 Let $A=R+Q$ be the canonical decomposition of $A$ with respect
 to $\varphi_T$.  Then  the operator $R:=[\text{\rm SOT}-\lim\limits_{k\to \infty} \varphi_T^k(A)]$
  is   selfadjoint  and $\varphi_T(R)=R$, and $Q$ is a pure solution of 
 the inequality $\varphi_T(X)\leq X$.
 Applying Corollary \ref{brat} to the operator $R$, we find 
 $B\in B(\cK)$ 
 in the commutant 
 of $\{V_i, V_i^*\}_{i=1}^\infty$
 such that $R=P_\cH B|\cH$. Using Theorem \ref{pure2}, we infer that $Q= KK^*$,
 as required. The proof of the converse is the same as the one from 
 Theorem \ref{main1}.
 \end{proof}

 The following result  provides a new insight and  an alternative proof of
  Corollary \ref{cons1}, which is based on the  Poisson transforms 
  of Section \ref{Poisson} and not on
   the noncommutative
   commutant lifting theorem.

\begin{proposition}\label{equal}
Let $\varphi_T$ be a w$^*$-continuous completely positive
linear map with $\varphi_T(I)\leq I$ and let $\varphi_V$ be its
 minimal isometric dilation. If $R\in B(\cH)$, $R\geq 0$, 
 is a solution of the equation $\varphi_T(X)=X$, then there is $A\in B(\cK)$ in the commutant of
  $\{V_i, V_i^*\}_{i=1}^\infty$ such that  $R=P_\cH A|\cH$.
  Conversely, if $A\in B(\cK)$ is such that $\varphi_V(A)=A$ then
  $R:=P_\cH A|\cH$  is a solution of the equation $\varphi_T(X)=X$
  
  Moreover, if  $\varphi_T(I)=I$, then the result remains true 
   if  $R\in B(\cH)$ is a selfadjoint  
      solution of the equation $\varphi_T(X)=X$.
\end{proposition}
\begin{proof}
Let $R\in B(\cH)$ be such that $0\leq R\leq I$ and  $\varphi_T(R)=R$.
According to Theorem \ref {poisson} there is a unique completely positive 
linear map
 $P_{\varphi, R}:C^*(S_1,\ldots, S_n)\to B(\cH)$
  such that
$P_{\varphi, R}(S_\alpha S_\beta^*) =T_\alpha RT_\beta^*$, \ $\alpha, 
\beta\in \FF_n^+$, and $P_{\varphi, R}(I)= R$.
Since $\varphi(I-R)=\varphi(I)-R\leq I-R$,
we can apply again Theorem \ref{poisson} and find a completely
 positive linear map
$P_{\varphi, I-R}:C^*(S_1,\ldots, S_n)\to B(\cH)$
  such that
$P_{\varphi, I-R}(S_\alpha S_\beta^*) =T_\alpha (I-R)T_\beta^*$, \ $\alpha, 
\beta\in \FF_n^+$, and $P_{\varphi, I-R}(I)= I-R$.
Hence, 
$P_{\varphi, I-R}=P_{\varphi, I}-P_{\varphi, R}$ is a completely
 positive linear map, where 
 $ P_{\varphi, I}$ is the Poisson transform associated with the row 
 contraction  $[T_1,\ldots, T_n]$. Therefore, 
 \begin{equation}\label{PP}
 0\leq P_{\varphi, R}\leq P_{\varphi, I}.
 \end{equation}
Notice that $P_{\varphi, I}(x)= P_\cH \pi(x) |\cH$, where $\pi$ 
is the representation of 
$C^*(S_1,\ldots, S_n)$ generated by the minimal isometric dilation 
 $[V_1,\ldots, V_n]$ on the Hilbert space $\cK\supseteq \cH$, i.e.,
$\pi(S_\alpha S_\beta^*)=V_\alpha V_\beta^*$.
On the other hand,
since $P_{\varphi, R}$ is a completely positive linear map, according to 
Stinespring theorem
\cite{S}, 
there is a representation
$\rho:C^*(S_1,\ldots, S_n)\to B(\cG)$ and an operator $W\in B(\cH, \cG)$ such that
$P_{\varphi, R}(x)= W^* \rho(x) W$, \ $x\in C^*(S_1,\ldots, S_n)$, and 
$\bigvee_{\alpha,\beta\in \FF_n^+} \rho(S_\alpha S_\beta^*) W\cH= \cG$.
Define the operator $C\in B(\cK, \cG)$ be setting
\begin{equation*}
C\left(\sum_{i=1}^k V_{\alpha_i} V_{\beta_i}^* h_{\alpha_i \beta_i}\right):=
\sum_{i=1}^k \rho(S_{\alpha_i} S_{\beta_i}^*) W h_{\alpha_i \beta_i},
\end{equation*}
where $\alpha_i, \beta_i \in \FF_n^+$, $h_{\alpha_i \beta_i}\in \cH$, 
and $k\in \NN$.
Using \eqref{PP}, we have
\begin{equation*}\begin{split}
\|\sum_{i=1}^k \rho(S_{\alpha_i} S_{\beta_i}^*) W h_{\alpha_i \beta_i}\|^2
&=
\sum_{i,j=1}^k \left<
P_{\varphi, R}(S_{\beta_j} S_{\alpha_j}^* S_{\alpha_i}
 S_{\beta_i}^*)h_{\alpha_i \beta_i}, h_{\alpha_j \beta_j}
\right>\\
&\leq \sum_{i,j=1}^k \left<
P_{\varphi, I}(S_{\beta_j} S_{\alpha_j}^* S_{\alpha_i}
 S_{\beta_i}^*)h_{\alpha_i \beta_i}, h_{\alpha_j \beta_j}\right>\\
 &=
 \|\sum_{i=1}^k V_{\alpha_i} V_{\beta_i}^* h_{\alpha_i \beta_i}\|^2.
 \end{split}
\end{equation*}
This shows that $C$ is well-defined and can be extended
 to a contraction from $\cK$ to $\cG$ with the properties
 $Ch =Wh, h\in \cH$ and 
 \begin{equation}\label{CVV}
 CV_{\alpha} V_{\beta}^*h= \rho(S_\alpha S_\beta^*)Wh,\quad 
 \alpha,\beta\in \FF_n^+, h\in \cH.
 \end{equation}
Hence, we deduce that 
\begin{equation}\label{CVRO}
 CV_{\alpha} V_{\beta}^*= \rho(S_\alpha S_\beta^*)C,\quad 
 \alpha,\beta\in \FF_n^+
\end{equation}
Indeed,
using \eqref{CVV}, we have
\begin{equation*}
CV_{\alpha} V_{\beta}^* V_\gamma h= \rho(S_\alpha S_\beta^* S_\gamma)Wh
=\rho(S_\alpha S_\beta^* )C V_\gamma h.
\end{equation*}
Now, setting  $B:= C^* C$, we infer that $0\leq B\leq I$ and
\begin{equation*}
\begin{split}
BV_{\alpha} V_{\beta}^*&=C^* C V_{\alpha} V_{\beta}^* 
= C^*\rho(S_\alpha S_\beta^*) C\\
&= V_{\alpha} V_{\beta}^*C^* C=V_{\alpha} V_{\beta}^*B.
\end{split}
\end{equation*}
Therefore, $B$ is in the commutant of $\{V_i,V_i^*\}_{i=1}^n$.
On the other hand, notice that 
\begin{equation}\begin{split}
\left< P_\cH BV_{\alpha} V_{\beta}^* h, h'\right>
&= \left< CV_{\alpha} V_{\beta}^*h, Ch'\right>\\
&=\left< \rho(S_\alpha S_\beta^*) Ch, Ch'\right>\\
&=\left< W^*\rho(S_\alpha S_\beta^*) Wh, h'\right>\\
&= \left< P_{\varphi, R}(S_\alpha S_\beta^*) h, h'\right>,
\end{split}
\end{equation}
for any $h,h'\in \cH$.
Therefore, $P_{\varphi, R}(S_\alpha S_\beta^*)= 
P_\cH BV_{\alpha} V_{\beta}^*|\cH$, \ 
$\alpha, \beta\in \FF_n^+$.
Hence, we have $R= P_{\varphi, R}(I)= P_\cH B|\cH$.
Since $B$ commutes with each $V_i$ and $V_i^*$, we have
$$
\varphi_V(B)=B^{1/2} \varphi_V(I) B^{1/2}\leq B
$$
if  $\varphi_V(I)\leq I$ , and $\varphi_V(B)=B$  if $\varphi_V(I)=I$.
We recall that $\varphi_V(I)=I$ if and only if  $\varphi(I)=I$.

For the converse, assume $A\in B(\cK)$ satisfies 
$\varphi_V(A)=A$, and let $R:= P_\cH A|\cH$.
Taking into accout that $[V_1,\ldots, V_n]$ is the minimal isometric dilation of 
$[T_1, \ldots, T_n]$, 
 we have
\begin{equation*}
\begin{split}
\left< \varphi(R)h, h'\right>&=\sum_{i=1}^n \left< AV_i^*h, V_i^* h'\right>=
\left< \varphi_V(A)h, h'\right>\\
&=\left< Ah, h'\right>=\left< Rh, h'\right>
\end{split}
\end{equation*}
for any $h,h'\in \cH$.

When $\varphi(I)=I$, the result of the theorem remains true if $R$ is 
a selfadjoint
operator  satisfying $\varphi(R)=R$. It is enough to apply the first
 part of the theorem to the positive operators $R_1, R_2$ in the Jordan
  decomposition $R=R_1-R_2$. Notice that, since $\varphi(I)=I$, 
   we have $\varphi(R_j)=R_j$, $j=1,2$.
  The proof is complete.
\end{proof}

\bigskip 
\section{ Common invariant subspaces for $n$-tuples of operators }
\label{Invariant}

We show that there is a strong connection between 
  the positive solutions of the operator inequality 
  $\varphi(X)\leq X$, where $\varphi$ is 
  a $w^*$-continuous completely positive linear map on $B(\cH)$ 
  defined by
  $\varphi(X):= \sum_{i=1}^n A_i X A_i^*, \quad X\in B(\cH)$,
  and the common invariant subspaces for the $n$-tuple of
   operators $\{A_i\}_{i=1}^n$.
   In this direction, we obtain invariant subspace theorems (eg.\,Theorem 
   \ref{fix}) and  Wold type  decompositions   for  
   $w^*$-continuous completely  positive linear maps on $B(\cH)$   
   (eg.\,Theorem \ref{wold1}).
  The latter results generalize 
    the classical
    Wold decomposition for isometries, as well as the one obtained
     in \cite{Po-isometric}
    for isometries with orthogonal ranges.
As in the previous  sections  we consider $n\in \NN$
 or $n=\infty$.
\begin{lemma}\label{inv-ker}
Let $\varphi_A$ be a $w^*$-continuous completely positive 
linear map on $B(\cH)$ given by $\varphi_A(X):= \sum_{i=1}^n A_iXA_i^*$. 
If $X\geq 0$ and $\varphi_A(X)\leq X$, then the subspace $\ker X$ is 
invariant under each $A_i^*$, \ 
$i=1,\ldots, n$. In particular, if $\cM$ is a subspace of $\cH$ and 
$\varphi_A(P_\cM)\leq P_\cM$, then $\cM$ is invariant under each $A_i$, \ 
$i=1,\ldots, n$. If $\varphi_A(I)\leq I$ and $\cM$ is reducing under each $A_i$, \ 
$i=1,\ldots, n$, then  $\varphi_A(P_\cM)\leq P_\cM$.
\end{lemma}

\begin{proof}
Since $X\geq 0$ and $\varphi_A(X)\leq X$, for any $h\in \ker X$, we have
$$
0\leq \sum_{i=1}^n \langle A_iXA_i^* h, h\rangle \leq \langle Xh,h\rangle=0.
$$
Hence $\|X^{1/2} A_i^* h\|=0$, \ $i=1,\ldots, n$, whence
$A_i^* h\in \ker X$.
As a particular case, 
if $\cM$ is a subspace of $\cH$ and 
$\varphi_A(P_\cM)\leq P_\cM$, then $\cM$ is invariant under each $A_i$, \ 
$i=1,\ldots, n$.
On the other hand, if 
 $\varphi_A(I)\leq I$ and $\cM$ is reducing under each $A_i$, \ 
$i=1,\ldots, n$, then $P_\cM A_i =A_i P_\cM$ and therefore
$$
\varphi_A(P_\cM)= P_\cM \varphi_A(I) P_\cM\leq P_\cM.
$$
The proof is complete.
\end{proof}

 \begin{corollary}\label{inv}
 Let $\varphi_T$ be a $w^*$-continuous completely positive 
linear map on $B(\cH)$ given by $\varphi_T(X):= \sum_{i=1}^n T_iXT_i^*$
 with $\varphi_T(I)\leq I$. If $ X\in B(\cH)$  is a positive operator 
 such that $\|X\|=1$ and 
 $\varphi_T(X)\geq X$ 
then the fixed-point set $\{h\in \cH: Xh=h\}$ is 
invariant under each $T_i^*$, \ 
$i=1,\ldots, n$.
 \end{corollary}
 \begin{proof}
 Notice that $I-X\geq 0$ and 
 $$
 \varphi_T(I-X)= \varphi_T(I)-\varphi_T(X)\leq I-X.
 $$
 Applying Lemma \ref{inv-ker}, to the positive operator 
 $I-X$, the result follows.
 \end{proof}

Let $\varphi$ be a $w^*$-continuous completely positive 
linear map on $B(\cH)$. The orbit of $X\in B(\cH)$ under  the semigroup 
generated by $\varphi$ is
 the sequence
$$
 \varphi^0(X), \varphi^1(X), \varphi^2(X), \ldots,
 $$ 
where $\varphi^0(X):=X$.
We say that the orbit of $X\in B(\cH)$ under $\varphi$ has a fixed point
 if there 
is $h\in \cH$ such that
$$
\varphi^0(X)h=\varphi^k(X)h, \quad k=1,2,\ldots.
$$ 

\begin{theorem}\label{fix}
Let $\varphi_A$ be a $w^*$-continuous completely positive 
linear map on $B(\cH)$
given by $\varphi_A(X):= \sum_{i=1}^n A_iXA_i^*$.
  Assume that   there is a positive operator 
$X\in B(\cH)$, $X\neq 0$,  such that $\varphi_A(X)\leq X$.
If one of the following statements holds, then 
there is a nontrivial invariant subspace under each $A_i$, $i=1,\ldots, n$:
\begin{enumerate}
\item[(i)]
$X$ is not injective;
\item[(ii)]
$X$ is not pure with respect to $\varphi_A$ and there is 
$h\in \cH$, $h\neq 0$, such that  $\lim\limits_{k\to\infty} \varphi^k_A(X)h=0$;
\item[(iii)] $\varphi_T(X)\neq X$, 
 and the orbit of 
$X$ under $\varphi_A$ has a nonzero fixed point. 
\end{enumerate}
\end{theorem}
\begin{proof} Due to Lemma \ref{inv-ker}, if $X$  is not injective, then
$(\ker X)^\perp$ is  invariant under each  each $A_i$, $i=1,\ldots, n$.
Now  assume that (ii) or (iii) holds. Let 
  $X= B+C$ be the canonical decomposition of $X$
 with respect to 
$\varphi_A$.
 According to Theorem \ref{decomp}, we have 
$$
B=\text{\rm SOT}-\lim_{k\to\infty} \varphi_A^k(X),
\quad \varphi_A(B)=B,
$$
 and $C$ is a pure solution of the inequality 
    $\varphi_A(X)\leq X$.
By Lemma \ref {inv-ker}, the subspaces $(\ker B)^\perp$ and $(\ker C)^\perp$ are
 invariant under each  each $A_i$, $i=1,\ldots, n$.
 On the other hand, we have
 $$
 \ker B= \{h\in \cH: \  \text{\rm SOT}-\lim_{k\to\infty} \varphi_A^k(X)h=0\}
 $$
 and 
 $$
 \ker C=\{ h\in \cH:\ \text{\rm SOT}-\lim_{k\to\infty} \varphi_A^k(X)h= Xh\}.
 $$
Since $\varphi_A(X)\leq X$, it is easy to see that 
$$
 \ker C=\{ h\in \cH:\  \varphi_A^k(X)h= Xh, \ k\in \NN\}.
 $$
Now, notice that    the condition (ii) (resp.~(iii)) holds, 
if and only if 
  the
the subspace $\ker B$ (resp.~$\ker C$) is nontrivial, and the result follows.
\end{proof}

 \smallskip
 Another consequence of Lemma \ref{inv-ker} is the following result which 
 was also obtained
  in \cite{BJKW}.
 \begin{corollary}\label{inv2}
 Let $\varphi_T$ be a $w^*$-continuous completely positive 
linear map on $B(\cH)$  such that $\varphi_T(I)= I$. Then $\cM\subseteq \cH$ is 
is an invariant subspace  under each \ $T_i$, \ 
$i=1,2,\ldots$,  if and only if 
$\varphi_T(P_\cM)\leq P_\cM$. Moreover, $\cM$ is reducing 
 under each $T_i$, \ 
$i=1,2,\ldots$, if and only if 
$\varphi_T(P_\cM)= P_\cM$.
 \end{corollary}

\begin{proof}  According to Lemma
\ref{inv-ker},
if $\varphi_T(P_\cM)\leq P_\cM$, then the subspace $\cH\ominus \cM$ 
is invariant under
under each $T_i^*$, \ 
$i=1,2,\ldots$. Conversely, assume 
$\cM$ is invariant under
under each $T_i$, \ 
$i=1,2,\ldots$.
Then $QT_i Q=QT_i$, where $Q:=I-P_\cM$. It is easy to see that
$$
\varphi_T(Q)Q=\varphi_T(I)Q=Q=Q\varphi_T(I)Q=Q\varphi_T(Q)Q.
$$
Since the operators $\varphi_T(Q)$ and $I-Q$ are positive and commuting, we have
$$
\varphi_T(Q)-Q=\varphi_T(Q)(I-Q)\geq 0.
$$
Hence, $\varphi_T(Q)\geq Q$. Since $\varphi_T(I)=I$, we infer that 
$\varphi_T(P_\cM)\leq P_\cM$. The last part of the corollary can be 
proved in a similar manner.
\end{proof}
\bigskip

 The following two propositions are needed to prove our Wold
  type decomposition theorem for $w^*$-continuous completely 
  positive linear maps on $B(\cH)$.


\begin{proposition}\label{phi}
Let $\varphi$ be a $w^*$-continuous positive 
linear map on $B(\cH)$ with $\|\varphi\|\leq 1$.
Then 
 $$
 \varphi^\infty(I):= \text{\rm SOT}-\lim_{k\to\infty} \varphi^k(I)
 $$
 exists and has the following properties:
 \begin{enumerate}
 \item[(i)]
 $0\leq \varphi^\infty(I)\leq I;$
 \item[(ii)]
 $\varphi(\varphi^\infty(I))=\varphi^\infty(I);$
 \item[(iii)]
 If $\varphi^\infty(I)\neq 0$, then $\|\varphi^\infty(I)\|=1;$
 \item[(iv)]
 If $\varphi^\infty(I)h\neq 0$, then $\varphi^k(I)h\neq 0$ for any $k\in \NN$.
 \end{enumerate}
\end{proposition}
\begin{proof}
The first two statements are particular cases of Theorem \ref{decomp}.
 We prove (iii). For any $h\in \cH$, 
$k\in \NN$, we have
$$
\left< \varphi^\infty(I)h,h\right>=\langle \varphi^k(\varphi^\infty(I))h,h\rangle
\leq \|\varphi^\infty(I)\| \langle \varphi^k (I)h,h\rangle.
$$
Taking the limit as $k\to\infty$, we obtain
$$
\left< \varphi^\infty(I)h,h\right>\leq 
\|\varphi^\infty(I)\|^2 \left< h,h\right>.
$$
Hence, $\|\varphi^\infty(I)^{1/2}\|\leq \|\varphi^\infty(I)^{1/2}\|^2$.
Since $\|\varphi^\infty(I)^{1/2}\|^2= \|\varphi^\infty(I)\|\leq 1$, 
we deduce that
$\|\varphi^\infty(I)^{1/2}\|=0$ or $\|\varphi^\infty(I)^{1/2}\|=1$, 
which proves (iii).
For the proof of (iv),  let  $h\in\cH$ be such that 
$\varphi^\infty(I)h\neq 0$ and assume that there is $k_0\in \NN$ with
$\varphi^{k_0}(I)h\neq 0$.
Since 
$$
0\leq \langle \varphi^{m+k_0}(I) h,h\rangle\leq 
\langle \varphi^{k_0}(I) h,h\rangle
=0,
$$
for any $m\in\NN$, we infer that  $\varphi^{m+k_0}(I) h=0$, $m\in\NN$. Hence 
$\varphi^\infty(I)h=0$, which is a contradiction. Therefore, 
 $\varphi^k(I)h\neq 0$ for any $k\in \NN$.
\end{proof}

Let us remark that if $\varphi(I)\leq I$, then
\begin{equation}\label{ker1}
\ker \varphi^\infty(I)=\{ h\in \cH:\ \lim_{k\to\infty} \varphi^k(I)h=0\}
\end{equation}
and
\begin{equation}\label{ker2}
\ker (I-\varphi^\infty(I))=\{ h\in \cH:\  \varphi^k(I)h 
=h, ~k\in \NN\}.
\end{equation}
\begin{proposition}\label{wold}
Let $\varphi$ be a   positive 
linear map on $B(\cH)$ with $\|\varphi\|\leq 1$.
Then 
$\cH$ admits a decomposition of the form
\begin{equation}
\label{wequ}
\cH=\cM\oplus \ker (I-\varphi^\infty(I))\oplus \ker \varphi^\infty(I),
\end{equation}
and $\cM=\{0\}$ if and only if $\varphi^\infty(I)$ is an orthogonal projection.
\end{proposition}
\begin{proof}
Since
$$
\cH= \overline{\text{\rm range}~\varphi^\infty(I)}\oplus \ker \varphi^\infty(I)
\ 
\text{ and } \ \ker (I-\varphi^\infty(I))\subseteq  
\text{\rm range}~\varphi^\infty(I),
$$
 we obtain relation \eqref{wequ}.
 On the other hand, since $\varphi^\infty(I)$ is a positive operator, 
 one can prove that
 $$
 \ker [\varphi^\infty(I)-  \varphi^\infty(I)^2]=\ker \varphi^\infty(I)
 \oplus \ker (I-\varphi^\infty(I)).
 $$
 On the other hand, since $ \ker [\varphi^\infty(I)-  \varphi^\infty(I)^2]=\cH$
 if and only if 
  $\varphi^\infty(I)$ is an orthogonal projection, 
  the result follows.
\end{proof}

 Now we can obtain the following Wold type decomposition.   
\begin{theorem}\label{wold1}
Let 
$\varphi_A$ be a $w^*$-continuous completely positive 
linear map on $B(\cH)$ given by $\varphi_A(X):= \sum_{i=1}^\infty A_iXA_i^*$
such that $\|\varphi_A\|\leq 1$. Then 
$\cH$ admits a decomposition of the form
\begin{equation*}
\cH=\cM\oplus \ker (I-\varphi_A^\infty(I))\oplus \ker \varphi_A^\infty(I),
\end{equation*}
and 
the subspaces $\ker (I-\varphi_A^\infty(I))$ and $\ker \varphi_A^\infty(I)$
 are invariant under each $A^*_i$, \ $i=1,\ldots, n$.
If, in addition,
$\varphi_A^\infty(I)$ is an orthogonal
projection,  then we have 
\begin{equation}\label{newwold}
\cH= \ker (I-\varphi_A^\infty(I))\oplus \ker \varphi_A^\infty(I),
\end{equation}
and 
the subspaces $\ker (I-\varphi_A^\infty(I))$ and $\ker \varphi_A^\infty(I)$
 are reducing for each $A_i$, \ $i=1,\ldots, n$.
 \end{theorem}
 \begin{proof} According to Proposition \ref{phi},  we have
 $\varphi_A(\varphi_A^\infty(I))=\varphi_A^\infty(I)$.
 Using  Lemma \ref{inv-ker} and Corollary \ref{inv}, we infer that 
 the subspaces $\ker \varphi_A^\infty(I)$ and  $\ker (I-\varphi_A^\infty(I))$
 are invariant under each $A^*_i$, \ $i=1,\ldots, n$. The rest
  of the proof follows from
  Proposition \ref{wold}. 
 \end{proof}
 Notice that,  in particular,  if   $A_i$ are isometries with orthogonal ranges, then 
 the decomposition \eqref{newwold} coincides with 
   the noncommutative
 Wold decomposition from \cite{Po-isometric}.

\bigskip
 
\section{ Similarity  of positive linear maps }\label{Similarity}

The main objectives of this section  are to 
 provide necessary and sufficient conditions 
for a $w^*$-continuous positive linear map $\varphi$ on $B(\cH)$ to be
similar to a positive  linear map $\lambda$ on $B(\cH)$ satisfying 
one of the following
properties:
\begin{enumerate}
\item[(i)]
$\lambda(I)=I$ (see Theorem \ref{simi});
\item[(ii)] 
$\|\lambda\|<1 $ (see Theorem \ref{simi2});
\item[(iii)] $\lambda$ is a pure completely positive linear map with
$\|\lambda\|\leq 1 $ (see Theorem \ref{sim-pure});
\item[(iv)]
$\lambda$ is a  completely positive linear map with
$\|\lambda\|\leq 1 $ (see Theorem \ref{simi4}).
\end{enumerate}
We show that these similarities are strongly related to the existence of 
invertible positive solutions of the operator inequality
$\varphi(X)\leq X$ or equation  $\varphi(X)= X$.

We say that two linear maps $\varphi, \lambda:B(\cH)\to B(\cH)$ are similar
if there is an invertible operator $R\in B(\cH)$ such that
\begin{equation}\label{SXS}
\varphi(RXR^*)= R\lambda(X)R^*, \quad \text{ for any } X\in B(\cH).
\end{equation}
Notice that relation \eqref{SXS} is equivalent to
\begin{equation}\label{SXS2}
\varphi = \psi_R\circ\lambda \circ \psi_R^{-1},
\end{equation}
where $\psi_R(X):= RXR^*$, \ $X\in B(\cH)$.

\begin{theorem} \label{simi}
Let $\varphi$ be a $w^*$-continuous  positive 
linear map on $B(\cH)$. Then the following statements are equivalent:
\begin{enumerate}
\item[(i)]
$\varphi$ if similar to a $w^*$-continuous  positive linear map 
$\lambda$ on $B(\cH)$,  with $\lambda(I)=I$.
\item[(ii)]
There exist  
 positive constants
 $0<a\leq b$ such that 
 \begin{equation*}
 aI \leq {\frac {\varphi^0(I)+\varphi^1(I)+\cdots + 
 \varphi^{k-1}(I)} {k}}\leq bI, \quad k\in \NN;
 \end{equation*}
\item[(iii)]
There exist  
 positive constants
 $0<a\leq b$ and an invertible positive operator $P\in B(\cH)$ such that 
 \begin{equation*}
 aI \leq {\frac {\varphi^0(P)+\varphi^1(P)+\cdots + 
 \varphi^{k-1}(P)} {k}}\leq bI, \quad k\in \NN;
 \end{equation*}
\item[(iv)]
There exist  
 positive constants
 $0<c\leq d$ such that
 $$
 cI\leq \varphi^k(I)\leq dI, \quad k\in \NN;
 $$
\item[(v)]
There exist  
 positive constants
 $0<c\leq d$  and an invertible positive operator $R\in B(\cH)$ 
 such that such that
 $$
 cI\leq \varphi^k(R)\leq dI, \quad k\in \NN;
 $$
 \item[(vi)]
 There exists  an invertible positive operator $Q\in B(\cH)$ 
 such that such that 
 $\varphi(Q)=Q  $.
\end{enumerate}
Moreover,    the operator $Q$ can be chosen such that $aI\leq Q\leq bI$.
\end{theorem}
\begin{proof}
First we prove that (i) $\Leftrightarrow$ (vi). Assume (i) holds, i.e.,   
$\varphi(RXR^*)= R\lambda(X)R^*$, where 
$\lambda(I)= I$ and
$R\in B(\cH)$ is an invertible operator
such that 
$aI\leq RR^*\leq bI$, for some constants $0<a\leq b$.
 Setting $X:=I$ and $Q:= RR^*$, 
we obtain $\varphi(Q)\leq Q$. Conversely, assume (vi) holds and define
$$
\lambda(X):= Q^{-1/2} \varphi(Q^{1/2} X Q^{1/2}) Q^{-1/2}, \quad X\in B(\cH).
$$
It is clear that $\lambda$ is a $w^*$-continuous positive linear map with $\lambda(I)=I$.
Moreover, we have 
$\varphi=\psi_R\circ \lambda\circ \psi_R^{-1}$, where $R:=Q^{1/2}$.

The implications (iv) $\Rightarrow$ (v) $\Rightarrow $ (iii) and  
(iv) $\Rightarrow$
(ii) $\Rightarrow$ (iii) are obvious. We prove that (i) $\Rightarrow$ (iv).
Assume 
$\varphi=\psi_R\circ \lambda\circ \psi_R^{-1}$, with $R\in B(\cH)$ invertible, and 
$\lambda(I)=I$.
Since all the maps are positive, we have
\begin{equation}\label{right}\begin{split}
\varphi^k(I)&\leq \|\psi_R^{-1} (I)\| \,\psi_R(\lambda^k(I))\\
&\leq 
\|\psi_R^{-1} (I)\| \| \psi_R(I)\|\, I =\|R^{-1}\|^2 \|R\|^2 \, I.
\end{split}
\end{equation}
On the other hand, we have 
\begin{equation*}
I= (\psi_R^{-1}\circ \varphi^k\circ \psi_R) (I)
\leq \|\psi_R(I)\| \, R^{-1} \varphi^k(I) {R^*}^{-1}, \quad k\in \NN.
\end{equation*}
Hence,  we obtain $RR^*\leq \|\psi_R(I)\| \, \varphi^k(I)$,  which implies 
\begin{equation}\label{left}
\varphi^k(I)\geq  \frac {1} {\|\psi_R(I)\|} RR^*\geq \frac{1} 
{\|R\|^2 \|R^{-1}\|^2} \, I.
\end{equation}
Putting  together relations \eqref{right} and \eqref{left}, we deduce
$$
\frac{1} 
{\|R\|^2 \|R^{-1}\|^2} \, I\leq \varphi^k(I)\leq \|R^{-1}\|^2 \|R\|^2 \,  I,
 \quad k\in \NN,
$$
which proves  (iv).

It remains to show that (iii) $\Rightarrow $ (vi). 
Let $P$ be an invertible positive operator such that (iii) holds. Since 
$\varphi$ is positive, we have
\begin{equation}\label{jplus}
\frac {1} {j+1} \varphi^j (P)\leq \frac {1} {j+1} \sum_{q=0}^j 
\varphi^q(P)\leq b I, \quad j\in \NN.
\end{equation}
On the other hand,  it is clear that, for any $j\leq k$,
$$
\varphi^k(P)=
(\varphi^{k-j}\circ \varphi^j)(P)\leq \|\varphi^j(P)\| \varphi^{k-j}(I).
$$
Hence, and using \eqref{jplus}, we infer that
\begin{equation} \label{bb}
\begin{split}
\varphi^k(P) \sum_{j=0}^k \frac {1} {j+1} &\leq
\sum_{j=0}^k \frac {1} {j+1}\|\varphi^{j}\| \, \varphi^{k-j}(I)\\
&\leq
b \sum_{j=0}^k  \varphi^{j}(I).
\end{split}
\end{equation}
Since $P$ is an invertible positive operator and $\varphi^j$ is positive,
 we have 
$I\leq \|P^{-1}\| P$ and 
$$
\varphi^j(I)\leq \|P^{-1}\| \varphi^j(P).
$$
Hence, and using again  the inequalities in (iii), we get
\begin{equation*} \begin{split}
\sum_{j=0}^k  \varphi^{j}(I)&\leq (k+1) \|P^{-1}\| \left( \frac {1} {k+1}
\sum_{j=0}^k \varphi^j(P)\right)\\
&\leq b(k+1)  \|P^{-1}\|\, I.
\end{split}
\end{equation*}
These inequalities together with \eqref{bb} imply 
$$
\left( \frac {1} {k} \varphi^k(P)\right) \sum_{j=0}^ k \frac {1} {j+1}
\leq \frac {b^2(k+1)} {k}  \|P^{-1}\|\, I,
$$
for any $k\in \NN$.
This implies that
$$
\sup_{k} \left\| \frac {1} {k} \varphi^k(P) 
\sum_{j=0}^ k \frac {1} {j+1}\right\| <\infty.
$$
   Hence,  we must have 
\begin{equation}\label{conv0}
\left\|\frac {1} {k} \varphi^k(P)\right\|\to 0, \quad \text{\rm as } k\to \infty.
\end{equation}
For each $k\geq 1$, we define the operator
$$
Q_k:=  \frac {1} {k} \sum_{j=0}^k \varphi^j(P).
$$
Since $\{Q_k\}_{k=1}^\infty$ is bounded sequence of operators and 
the closed unit ball of $B(\cH)$ is weakly compact, there is a subsequence 
 of $\{Q_k\}_{k=1}^\infty$
 weakly
 convergent   to an operator $Q\in B(\cH)$. 
 Due to (iii),  the operator $Q_k$ is  positive  and satisfies
 $aI\leq Q_k\leq bI$. Therefore, $Q$ is an invertible positive 
 operator satisfying the same inequalities.
 Since 
 $$Q_k-\varphi(Q_k)= \frac {1} {k} P- \frac {1} {k} \varphi^k(P)
 $$
 and taking ito account  \eqref{conv0}, we get $\|Q_k-\varphi(Q_k)\|\to 0$,
  as $k\to\infty$.
 Now, using the fact that the $w^*$ and the weak
  topologies coincide on bounded sets of $B(\cH)$, and that $\psi$
   is $w^*$-continuous,
  we infer that $\varphi(Q)= Q$.
  Therefore,  (vi) holds and the proof is complete.
\end{proof}

Let us consider an application of this theorem and the results
 from Section \ref{Poisson}.
Let $\{A_i\}_{i=1}^\infty\subset B(\cH)$ $(n\in\NN \text{ \rm or } n=\infty)$ 
be a sequence of operators such that 
$$
\varphi_A(X):= \sum_{i=1}^n  A_iXA_i^*, \quad X\in B(\cH),
$$
is a $w^*$-continuous completely positive linear map on $B(\cH)$.
 According to  Theorem \ref{simi},
 if there exist  some constants $0<a\leq b$ such that
 \begin{equation} \label{asumb}
 aI\leq \frac{1}{k} \sum_{j=0}^{k-1} \sum_{|\alpha|=j}
  A_\alpha A_\alpha^*\leq bI,\quad k\in \NN,
 \end{equation}
 then there is an invertible positive operator $Q\in B(\cH)$
  such that $\varphi_A(Q)=Q$
 and $aI\leq Q\leq  bI$.
   Therefore, 
 we can apply Theorem \ref {poisson} and  Theorem \ref{Poisson-Cuntz} 
 to the map $\varphi_A$.
 In particular one can deduce that
 $$
 \|p(A_1,\ldots, A_n) Q^{1/2}\|\leq \|Q^{1/2}\| \|p(S_1,\ldots, S_n)\|,
 $$
 for any polynomial $p(S_1,\ldots, S_n)$ in 
 the noncommutative disc algebra $\cA_n$.
  Hence, we obtain the inequality
 \begin{equation}\label{cp-von}
 \|p(A_1,\ldots, A_n)\|\leq \|Q^{1/2}\| \|Q^{-1/2}\| \|p(S_1,\ldots, S_n)\|,
 \end{equation}
 which can be extended to matrices over $\cA_n$.
 Therefore, if \eqref{asumb} holds, then 
  the  homomorphism 
 $\Phi:\cA_n\to B(\cH)$,
    defined by $\Phi(p):=
 p(A_1,\ldots, A_n)$,  is completely bounded and 
 $\|\Phi\|_{cb}\leq \sqrt{\frac {b} {a}}$.
 We remark that,  the inequality \eqref {cp-von} remains true if we 
 replace the left creation operators
 by a set of generators of the Cuntz algebra $\cO_n$.


\begin{corollary}\label{inve}
Let $\varphi$ be a w*-continuous positive linear map on $B(\cH)$ such that
\begin{equation}\label{so}
\varphi^\infty(I):= \text{\rm SOT}-\lim_{k\to \infty} \varphi^k(I)
\end{equation}
exists. Then $\varphi$ is similar to a positive linear map $\psi$ such that $\psi(I)=I$ 
if and only if $\varphi^\infty(I)$ is invertible.
\end{corollary}
\begin{proof}
Since  $\varphi$ is  $w^*$-continuous  and the limit \eqref{so} exists, we have
$\varphi(\varphi^\infty(I))=\varphi^\infty(I)$. If  $\varphi^\infty(I)$ is 
invertible, then 
the result follows from Theorem \ref{simi}. Conversely,  assume that $\varphi$
is similar to a positive linear map $\psi$ such that $\psi(I)=I$.
 Using again
Theorem \ref{simi}, we have $cI\leq \varphi^k(I)\leq dI$ for any $k\in\NN$.
Hence, $\varphi^\infty(I)$ is invertible. 
\end{proof}
Let us remark that the limit in \eqref{so} exists in 
the particular case when $\varphi(I)\leq I$.

In what follows, we find sufficient conditions 
 for the existence of an injective operator in  $C_=(\varphi)$.
We say that a   positive 
linear map $\varphi$ on $B(\cH)$ is power bounded if there is
 a  constant $M>0$
such that
$$
\|\varphi^k\|\leq M, \quad k\in \NN.
$$ 

\begin{lemma}\label{power1}
 Let $\varphi$ be a $w^*$-continuous  positive 
linear map on $B(\cH)$ such that $\varphi$ is power bounded and
$
\langle \varphi^k(I)h,h\rangle 
$
 does not converge to zero for any $h\in \cH$, $h\neq 0$.
 Then there is an injective positive solution of the   equation
 $\varphi(X)=X$.
\end{lemma}
\begin{proof}
First let us prove that if $Y\in B(\cH)$, $Y\geq 0$, then the set
\begin{equation}\label{conv}
\{X\geq0:\ \varphi(X)=X\} \cap 
\overline{\text{\rm conv}}^w \{\varphi^k(Y):\ k=0,1,\ldots\}
\end{equation}
is nonempty, where  $\overline{\text{\rm conv}}^w$ stands for the 
weakly closed convex hull.
Since $\|\varphi^k\|\leq M$, $k\in \NN$, the sequence of Cesaro means
$$
\sigma_k(Y):= {\frac {\varphi^0(Y)+\varphi^1(Y)+\cdots + 
 \varphi^{k-1}(Y)} {k}},\quad k=1,2,\ldots, 
 $$
 is bounded. Therefore, there is a subsequence
  $\{\sigma_{n_k}(Y)\}_{k=1}^\infty$ weakly
  convergent to an operator
   $Z\in 
  \overline{\text{\rm conv}}^w \{\varphi^k(Y):\ k=0,1,\ldots\}$.
Since $\varphi$ is $w^*$-continuous  positive 
linear map on $B(\cH)$ and 
$$
\|\varphi(\sigma_k(Y))-\sigma_k(Y)\|\leq 
\frac {\|Y\|} {k}+ \frac {\|\varphi^k(Y)\|} {k}\leq \frac {(M+1) \|Y\|} {k},
$$
for any $k=1,2,\ldots$, 
we infer that $\varphi(Z)=Z$. Hence the set \eqref{conv} is nonempty.
Now, set $Y:=I_\cH$ and let $Z$ be in the set \eqref{conv}.
Let $h\in \cH$, $h\neq 0$, and assume that 
$
\langle \varphi^k(I)h,h\rangle 
$
 does not converge to zero, as $k\to\infty$. This implies that there is a constant $C>0$
 such that 
 \begin{equation}\label{const}
 \langle \varphi^k(I)h,h\rangle \geq C,\quad \text{ for any } k=0,1,\ldots.
 \end{equation}
Indeed, if this were not true, then there would exist a subsequence $\{n_k\}$ 
such that 
$\langle \varphi^{n_k}(I)h,h\rangle \to 0$, as $k\to\infty$.
Notice that
\begin{equation*}\begin{split}
\langle \varphi^m(I)h,h\rangle 
&= 
\langle \varphi^{n_k}(\varphi^{m-n_k}(I))h,h\rangle \\
&\leq \|\varphi^{m-n_k}(I)\| \langle \varphi^{n_k}(I)h,h\rangle\\
&\leq M \langle \varphi^{n_k}(I)h,h\rangle.
\end{split}
\end{equation*}
Therefore, $\langle \varphi^{m}(I)h,h\rangle\to 0$ as $m\to\infty$,
 which is a contradiction.
Hence, relation  \eqref{const} holds and we have
$$
\langle \sum_{k\geq 0} \gamma_k\varphi^{k}(I)h,h\rangle\geq C
$$
for any finitely supported sequence $\{\gamma_n\}$ of positive numbers with
$\sum_{k\geq 0} \gamma_k=1$.  Consequently, $\langle Zh,h\rangle \geq C$, which 
shows that $Z$ is an injective operator. 
Moreover, we have 
$$\ker Z=\{h:\  \langle \varphi^k(I)h,h\rangle\to 0, \ \text{\rm as } k\to\infty\}.
$$
The proof is complete.
\end{proof}

 \begin{corollary}\label{findim} 
 If $\cH$ is finite dimensional and 
$\varphi$ is  a $w^*$-continuous positive 
linear map on $B(\cH)$ such that  $\varphi$ is power bounded and
$
\langle \varphi^k(I)h,h\rangle 
$
 does not converge to zero for any $h\in \cH$, $h\neq 0$, then
 $\varphi$ is similar to  a  positive 
linear map $\psi$ with $\psi(I)=I$.
\end{corollary}
\begin{proof}
According to Lemma \ref{power1}, there is an injective positive 
 operator $Z\in B(\cH)$ such that $\varphi(Z)=Z$. Since 
 $\cH$ is finite dimensional, $Z$ is invertible. Using Theorem \ref{simi},
 the result follows.
\end{proof}

Now, we present a few results concerning the similarity
 of positive linear maps with contractive (resp.~strictly contractive) ones.

\begin{lemma}\label{strict}
Let $\varphi$ be a positive linear map on $B(\cH)$.
Then $\varphi$ if similar to a    positive linear map 
$\psi$  with $\|\psi\|<1$  if and only if there is an 
invertible positive operator $R\in B(\cH)$ such that $R-\varphi(R)$
 is positive and invertible.
\end{lemma}
\begin{proof}
If  $Q$ is an invertible operator such that 
$\|\psi_Q^{-1}\circ \varphi \circ \psi_Q\|<1$, then 
we have 
$$
\|Q^{-1} \varphi(QQ^*) {Q^*}^{-1}\|<c,
$$
 for some positive constant $c<1$.
Hence, we get 
$\varphi(QQ^*)\leq c \,QQ^*$. Setting $R:= QQ^*$, we have
$$
R-\varphi(R)=(1-c) R,
$$
which is an invertible positive operator.
Conversely,
assume that $R\in B(\cH)$ is an invertible positive operator and 
\begin{equation}\label{RR}
R-\varphi(R)\geq b I
\end{equation}
 for some constant $b>0$.
 Let $a$ be such that $0<a<1$ and $a<\frac {b} {\|R\|}$.
 Since $R\leq \|R\|\, I\leq \frac {b} {a} \, I$,
 we infer that 
 \begin{equation}
 \label{RRRR}
 R-b I\leq R-a R.
 \end{equation}
Using  relations \eqref{RR} and \eqref{RRRR}, we obtain
$$
(1-a) R-\varphi(R)\geq  R-b I-\varphi(R)\geq 0.
$$
Hence, $\varphi(R)\leq (1-a) R$, which implies
$\|\psi_{R^{1/2}}^{-1}\circ \varphi \circ \psi_{R^{1/2}}\|<1$.
This completes the proof.
\end{proof}

\begin{proposition} \label{power}
Let $\varphi$ be a    positive 
linear map on $B(\cH)$. If there is $m\in \NN$ such $\varphi^m$ is
 similar to a positive linear map
$\psi$ with $\|\psi\|\leq 1$ (resp. $\|\psi\|< 1$), then
 $\varphi$   is similar to a positive linear map $\psi'$
 with $\|\psi'\|\leq 1$ (resp. $\|\psi'\|< 1$).
 \end{proposition}
 
 \begin{proof}
 If $\varphi^m$ is
 similar to a positive linear map
$\psi$ with $\|\psi\|\leq 1$, then there exists an 
invertible positive operator $Q$ such that $\varphi^m(Q)\leq Q$.
Denote
$$
P:= Q+\varphi(Q)+\cdots +\varphi^{m-1}(Q)
$$
and notice that $P$ is an invertible positive operator.
Moreover, we have
\begin{equation*}\begin{split}
\varphi(P)&=\varphi(Q)+\varphi^2(Q)+\cdots +\varphi^{m}(Q)\\
&\leq Q+\varphi(Q)+\cdots +\varphi^{m-1}(Q)=P.
\end{split}
\end{equation*}
Setting $\lambda(X):=P^{-1/2} \varphi(P^{1/2} X P^{1/2})P^{-1/2}$, 
the result follows.
Notice that if $\|\psi\|<1$, then, according to Lemma \ref{strict}, the operator
$Q-\varphi^m(Q)$ is invertible and positive.
Since $P-\varphi(P)= Q-\varphi^m(Q)$, we can apply again Lemma \ref{strict}
 to complete the proof.
 \end{proof}

\begin{proposition} \label{simi3}
Let $\varphi$ be a $w^*$-continuous  positive 
linear map on $B(\cH)$.
If there exist  
 positive constants
 $0<a\leq b$   and a  positive operator 
  $P\in B(\cH)$ such that 
 \begin{equation} \label{alfa}
 aI \leq  \sum_{k=0}^\infty \varphi^k(P) \leq bI,
 \end{equation}
 then 
 $\varphi$ if similar to a $w^*$-continuous  positive linear map 
$\psi$ on $B(\cH)$,  with $\|\psi\|\leq1$.
Moreover, if, in addition, $P$ is invertible, then 
$\varphi$ if similar to a $w^*$-continuous  positive linear map 
$\psi$  with $\|\psi\|<1$.
\end{proposition}
\begin{proof}
Setting $Q_m:= \sum_{k=0}^{m-1} \varphi^k(P)$, we have
\begin{equation}\label{QQ}
Q_{m+1}= Q_m+\varphi^m(P).
\end{equation}
Since $\{Q_m\}_{m=1}^\infty$  is a bounded monotone sequence 
of positive operators, it converges strongly to an operator $Q$.
Due to \eqref{alfa}, $Q$ is an invertible positive operator. 
According to \eqref{QQ},
$\varphi^m(P)\to 0$ strongly, as $m\to\infty$.
Since 
$$
Q_m-\varphi(Q_m)=P-\varphi^m(P),
$$
we get $Q-\varphi(Q)=P\geq 0$.
If $P$ is invertible, then we can apply Lemma \ref{strict} 
to complete the proof.
\end{proof}

Let $\varphi$ be a  positive 
linear map on $B(\cH)$. The spectral radius of $\varphi$ is defined by setting
$$
\text{\rm r}(\varphi):= \lim\limits_{k\to \infty} 
\|\varphi^k\|^{{\frac {1} {k}}}.
$$

\begin{lemma} \label{simil}
Let $\varphi$ be a    positive 
linear map on $B(\cH)$. Then the following statements are equivalent:
\begin{enumerate}
\item[(i)] $\text{\rm r}(\varphi)<1$;
\item[(ii)] $\lim\limits_{k\to \infty} \|\varphi^k\|=0$;
\item[(iii)] $\sum_{k=1}^\infty \|\varphi^k\|^p$ is convergent
 for any $p>0$.
\end{enumerate}
\end{lemma}
\begin{proof}
If (i) holds, then for any $a\in (r(\varphi), 1)$ there 
is $m\in \NN$ such that $\|\varphi^k\|\leq a^k$ for any $k\geq m$.
this clearly implies conditions (ii) and (iii).
Since 
$$
r(\varphi^n)=\lim_{k\to \infty}\left(\|\varphi^{nk}\|^{\frac {1}{nk}}\right)^n
=r(\varphi)^n
$$
 and $r(\varphi^n)\leq \|\varphi^n\|$ for any $n\in \NN$, it is clear that
(iii) implies (i).
\end{proof}

Now we can characterize those $w^*$-continuous positive
 linear maps on $B(\cH)$ which are similar 
to strictly contractive ones.

\begin{theorem} \label{simi2}
Let $\varphi$ be a $w^*$-continuous  positive 
linear map on $B(\cH)$. Then the following statements are equivalent:
\begin{enumerate}
\item[(i)]
$\varphi$ if similar to a $w^*$-continuous  positive linear map 
$\psi$ on $B(\cH)$,  with $\|\psi\|<1$.
\item[(ii)]
For any 
 invertible positive operator $R\in B(\cH)$   
  the equation 
  \begin{equation}\label{eq}
  X-\varphi(X)=R
  \end{equation}
  has an invertible positive  solution in $B(\cH)$.
 \item[(iii)] There esists an invertible positive operator $Q\in B(\cH)$
  such that $\varphi(Q)\leq Q$ and
 $Q-\varphi(Q)$ is invertible.
 \item[(v)] $\lim\limits_{k\to \infty} \|\varphi^k(I)\|=0$.
 \item[(vi)] $r(\varphi)<1$. 
\end{enumerate}
 Moreover, in this case the positive solution of  the equation \eqref{eq} is
  unique and given  by 
  $$
  X=\sum_{k=0}^\infty \varphi^k(R),
  $$
  where the convergence is in the uniform topology.
\end{theorem}
\begin{proof}
The equivalence  (i) $\Leftrightarrow$ (iii)
( resp. (v) $\Leftrightarrow$ (vi) ) was proved in Lemma \ref{strict}
(resp. Lemma \ref{simil}).
 In what follows we prove that 
 (ii) $\Rightarrow$ (i) $\Rightarrow$ (vi) 
  $\Rightarrow$ (ii).
  Assume (ii) holds. Let $Q\in B(\cH)$ be an invertible positive operator
   such that 
  $Q-\varphi(Q)=R$. Using Lemma \ref{strict}, we infer (i). Now, we assume (i).
  Then there is an invertible operator $Q$ such that 
  $\|\psi_Q^{-1}\circ \varphi \circ \psi_Q\|<1$.
  On the other hand,
  $$
  r(\varphi)=r(\psi_Q^{-1}\circ \varphi \circ \psi_Q)\leq 
  \|\psi_Q^{-1}\circ \varphi \circ \psi_Q\|<1.
  $$
  According to Lemma \ref{simil}, the latter condition is equivalent
   to the fact that the series  $\sum_{k=1}^\infty \|\varphi^k\| $
    is convergent. Then, for any invertible positive operator $R\in B(\cH)$,
    we have 
    $$
    \frac {1} {\|R^{-1}\|} \, I \leq R\leq \sum_{k=0}^\infty \varphi^k(R)
    \leq \left(\|R\| \sum_{k=0}^\infty \|\varphi^k\|\right)\, I.
    $$
   According to Proposition \ref{simi3} (see the proof),
   there is an invertible operator $Q$ such that $Q-\varphi(Q)=R$, so that 
   (ii) holds.
   
   To prove the last part of the theorem,
   let $X\geq 0$ be an invertible operator such 
   $X-\varphi(X)=R$, where $R\geq 0$ is a fixed  invertible operator.
   Let $X_k:= \sum_{j=0}^{k-1} \varphi^j(R)$, $k\in \NN$. Since 
   $$
   \varphi^j(R)=\varphi^j(X)-\varphi^{j+1}(X), \quad j=0,1,\ldots,
   $$
    we have
   $X_k=X-\varphi^k(X)$, for any $k\in \NN$. Since $\|\varphi^k\|\to 0$, 
   as $k\to \infty$, we have 
   $$
   0\leq \|X_k-X\|=\|\varphi^k(X)\|\leq \|X\| \|\varphi^k\|\to 0,
   $$
   as $k\to \infty$.
   Therefore, $X_k$ converges to $X$ uniformly, as $k\to \infty$,  and $X$ is the unique solution
   of the inequality $X-\varphi(X)=R$.
   The proof is complete.
   
\end{proof}

\begin{corollary} \label{radius}
Let $\varphi$ be a $w^*$-continuous  positive 
linear map on $B(\cH)$. Then
\begin{equation}\label{rad}
\text{\rm r}(\varphi):= \inf_{Q}\|\psi_Q\circ \varphi\circ \psi_Q^{-1}\|,
\end{equation}
where the infimum is taken over all invertible operators $Q\in B(\cH)$, and 
$\psi_Q(X):=QXQ^*$.
\end{corollary}
\begin{proof}
Let $\epsilon>0$ and denote $\varphi_\epsilon:= \frac {1} {r(\varphi)+\epsilon}
\cdot \varphi$.
Since 
$$
r(\varphi_\epsilon)=\frac {r(\varphi)} 
{r(\varphi)+\epsilon}<1,
$$
 we can apply Theorem \ref{simi2} to deduce 
that $\varphi_\epsilon$ 
is similar to a strictly contractive positive map on $B(\cH)$. Therefore, there
is an invertible operator $R_\epsilon\in B(\cH)$ such that 
$\|\psi_{R_\epsilon}\circ \varphi_\epsilon\circ \psi_{R_\epsilon}^{-1}\|<1$.
Hence 
\ $ \frac {1} {r(\varphi)+\epsilon}
\,\|\psi_{R_\epsilon}\circ \varphi\circ \psi_{R_\epsilon}^{-1}\|<1
$
and we can deduce that
$$
\inf_{Q}\|\psi_Q\circ \varphi\circ \psi_Q^{-1}\|\leq 
\|\psi_{R_\epsilon}\circ \varphi\circ \psi_{R_\epsilon}^{-1}\|
\leq r(\varphi)+\epsilon,
$$
for any $\epsilon>0$.
On the other hand, we have
$$
r(\varphi)= r(\psi_Q\circ \varphi\circ \psi_Q^{-1})
\leq \|\psi_Q\circ \varphi\circ \psi_Q^{-1}\|
$$
for all invertible operators $Q\in B(\cH)$. Now the equation \eqref{rad}
follows and the proof is complete.
\end{proof}

The next result provides necessary and sufficient conditions for a
$w^*$-continuous completely positive linear map to be  similar
 to one which is pure and contractive.
We recall that a positive 
linear map $\varphi$ on $B(\cH)$ is pure if $\varphi^k(I)\to 0$ strongly,
 as $k\to\infty$.

\begin{theorem}\label{sim-pure}
Let $\varphi$ be a $w^*$-continuous completely positive 
linear map on $B(\cH)$. The following statements are equivalent:
\begin{enumerate}
\item[(i)] $\varphi$ is similar to a pure completely positive linear
 map $\psi$ with $\|\psi\|\leq 1$;
 \item[(ii)] There exist two constants 
 $0<a\leq b$   and a  positive operator 
  $R\in B(\cH)$ such that 
 \begin{equation}\label{ab}
 aI \leq  \sum_{k=0}^\infty \varphi^k(R) \leq bI;
 \end{equation}
 \item[(iii)] There is an invertible  pure solution of the inequality
 $\varphi(X)\leq X$.
\end{enumerate}
\end{theorem}
\begin{proof} Let us prove that (ii) $\Rightarrow$ (i).
Assume that
(ii) holds.
 Since 
$\varphi$ is a $w^*$-continuous completely positive 
linear map on $B(\cH)$,
 there exists a sequence $\{A_i\}_{i=1}^n$ ($n\in\NN$ or $n=\infty$)
such that $\varphi (X)=\sum_{i=1}^nA_i XA_i^*$.
Let $W:\cH\to F^2(H_n)\otimes \cH$ be defined by
$$
Wh:= \sum_{k=0}^\infty \sum_{|\alpha|=k} e_\alpha\otimes 
R^{1/2} T_\alpha^*h,\quad, h\in \cH.
$$
Since $\|Wh\|^2=\sum_{k=0}^\infty \langle \varphi^k(R)h,h\rangle$ and 
relation \eqref{ab} holds, the range of $W$ is a closed subspace of 
$F^2(H_n)\otimes \cH$. Notice that 
$$
WA_i=(S_i^*\otimes I_\cH)W,\quad i=1,\ldots, n,
$$ and therefore the range of $W$  is invariant under each $S_i^*\otimes I_\cH$,
\ $i=1,\ldots, n$.
Since  the operator $W:\cH\to \text{\rm range}\, W$ is invertible, we have
\begin{equation}\label{WTW}
WA_i^*W^{-1}= (S_i^*\otimes I_\cH)|\text{\rm range}\, W,\quad i=1,\ldots, n. 
\end{equation}
Hence, $\varphi$ is similar to $\varphi_T$, where 
$\varphi_T (Y):=\sum_{i=1}^n T_i Y T_i^*$
and 
$$
T_i:= P_{ \text{\rm range}\, W} (S_i\otimes I_\cH) | 
 \text{\rm range}\, W,\quad i=1,\ldots, n.
$$
We clearly have  $\|\varphi_T\|\leq 1$ and 
$\varphi_T^k(I)\to 0$~ strongly, as $k\to \infty$. This proves (i).

Assume  now that $\varphi$ is similar to $\lambda$  such that
 $\|\lambda\|\leq 1$ and  $\lambda$ is pure. Hence, there is an 
 invertible operator $Q\in B(\cH)$ such that
 $\varphi=\psi_Q\circ \lambda \circ\psi_Q^{-1}$.
 Therefore, $D:= QQ^*$ satisfies $\varphi(D)\leq D$.
 Since $\lambda$ is pure and 
 $\psi_Q^{-1}\circ \varphi^k \circ \psi_Q=\lambda^k$, we deduce that
 $$
 \varphi^k(D)=Q\lambda^k(I) Q^*\to 0 \quad  \text{\rm strongly,  as } k\to \infty.
 $$
Now, let $R:= D-\varphi(D)$ and notice that $\sum_{k=0}^\infty \varphi^k(R)=D$,
which is invertible and positive. This implies  relation \eqref{ab}.

We already proved that (i) $\Rightarrow$ (iii). The implication 
(iii) $\Rightarrow$ (i) is easy. Indeed, asuume that $D\geq 0$ is an invertible pure solution of the inequality $\varphi(X)\leq X$. 
Then $\varphi(D)\leq D$ and, if we define
$$\lambda(X):= D^{-1/2} \varphi(D^{1/2} X D^{1/2}) D^{-1/2},
$$
we have $\lambda(I)\leq I$ and $\lambda^k(I)= 
D^{-1/2} \varphi(D) D^{-1/2}\to  0$  strongly, as $k\to\infty$.
The proof is complete.
\end{proof}

\begin{corollary} \text{\rm (\cite{Po-models})}
Let $\varphi$ be a $w^*$-continuous completely positive 
linear map on $B(\cH)$.
If there is a constant $b>0$ such that 
\begin{equation*}
 \sum_{k=0}^\infty \varphi^k(I) \leq bI,
 \end{equation*}
 then $\varphi$ is similar to a pure completely positive linear map $\lambda$
 with $\|\lambda\|\leq 1$.
 
\end{corollary}

Notice that if 
$\varphi$ is  a  positive 
linear map on $B(\cH)$ with $\|\varphi\|\leq 1$, then 
$$
I=\varphi^0(D)+\varphi^1(D)+\cdots + 
 \varphi^{k-1}(D)+ \varphi^k(I), \quad k\in \NN,
 $$
 where $D:= I-\varphi(I)$.
 In what follows, we show that a perturbation of this equality
 provides a characterization of those  
 $w^*$-continuous completely positive 
linear maps on $B(\cH)$ which are similar to contractive ones.

\begin{theorem} \label{simi4}
Let $\varphi_A$ be a $w^*$-continuous completely positive 
linear map on $B(\cH)$, given by $\varphi_A(X):=\sum_{i=1}^n A_iXA_i^*$ 
$(n\in \NN \ \text{\rm  or } n=\infty)$. 
Then the following statements are equivalent:
\begin{enumerate}
\item[(i)]
$\varphi_A$ if similar to a $w^*$-continuous completely  positive linear map 
$\psi$ on $B(\cH)$,  with $\|\psi\|\leq 1$;
\item[(ii)] There is an invertible positive operator $R\in B(\cH)$ such that
$\varphi_A(R)\leq R$;
\item[(iii)]
There exist  
 positive constants
 $0<a\leq b$   and a  positive operator 
  $D\in B(\cH)$ such that 
 \begin{equation*}
 aI \leq  \varphi_A^0(D)+\varphi_A^1(D)+\cdots + 
 \varphi_A^{k-1}(D)+ \varphi_A^k(I) \ \leq bI, \quad k\in \NN;
 \end{equation*}
\item[(iv)] The map $\Phi:\cA_n\to B(\cH)$ defined by $\Phi(p):=
 p(A_1,\ldots, A_n)$  is completely bounded, where $\cA_n$ is the noncommutative 
 disc algebra.
\end{enumerate}
  Moreover, if \text{\rm (iii)} holds, then ~$\|\Phi\|_{cb}\leq 
  \sqrt{ \frac {b} {a}}$.
\end{theorem}
\begin{proof}
The equivalence (i) $\Leftrightarrow$ (ii) is clear. Let us prove that 
(ii) $\Rightarrow$ (iii). Assume that $R$ is an invertible operator such that
$0\leq R\leq I$ and $\varphi_A(R)\leq A$.
Let  $D:= R-\varphi_A(R) $ and notice that
\begin{equation*} \begin{split}
\varphi_A^0(D)+\varphi_A^1(D)+\cdots + 
 \varphi_A^{k-1}(D)+ \varphi_A^k(I) &=
 R-\varphi^k_A(R) + \varphi_A^k(I)\\
 &\geq R\geq \frac {1} {\|R^{-1}\|}\, I.
\end{split}
\end{equation*}
On the other hand, since $\varphi_A$ is similar to a positive linear map of norm less than one,
there is a constant $M>0$ such that $\varphi_A^k(I)\leq M$, $k\in \NN$.
Moreover,
$$
R-\varphi^k_A(R) + \varphi_A^k(I) \leq R+\varphi^k_A(R)\leq  (1+M) I
$$
and  (iii) holds.
The implication (iii) $\Rightarrow$ (i) follows from Proposition 2.6 of
\cite {Po-models}. Indeed, if (iii) holds, then using \cite{Po-models},
 we find a sequence
 $\{T_i\}_{i=1}^n$ of operators such that $[T_1,\ldots, T_n]$ is a row
 contraction and $A_i= Y T_i Y^{-1}$, where $Y\in B(\cH)$ is an invertible
  positive operator and $\sqrt{a} I\leq Y\leq \sqrt{b} I$.
  Set $\varphi_T(X):= \sum_{i=1}^n T_i XT_i^*$ and notice that
   $\varphi_T(I)\leq I$ and $\varphi_A$ is similar to $\varphi_T$.
   On the other hand, we have 
   $$
   p(A_1,\ldots, A_n)= Yp(T_1,\ldots, T_n)Y^{-1}
   $$
   for any polynomial $p(S_1,\ldots, S_n)\in \cA_n$.
   Hence, we infer 
   (iv) and, using the noncommutative von Neumann inequality \cite{Po-von},
    we get
   $$
   \|\Phi\|_{cb}\leq\|Y\|\|Y^{-1}\|\leq \sqrt{\frac {b}{a}}.
   $$
   If we assume that (iv) holds, then using  Paulsen's result \cite{P} 
   and \cite{Po-disc}, we infer that
   $\{A_i\}_{i=1}^n$ is simultaneously similar to  $\{T_i\}_{i=1}^n$.
   This implies that $\varphi_A$ is similar to $\varphi_T$. 
   The proof is complete.
\end{proof}

\section{ Numerical invariants  for Hilbert modules over free 
semigroup algebras }
\label{Invariants}

 The Poisson transforms of Section \ref{Poisson} are used in  this section
    to define certain  numerical  invariants associated 
 with  (not necessarily contractive) Hilbert modules  
 over  
  the  free semigroup algebra  $\CC \FF_n^+$. 
  Any Hilbert module  $\cH$ over $\CC\FF_n^+$ corresponds to a unique 
  $w^*$-continuous
   completely positive map $\varphi$ on $B(\cH)$ and therefore
    to a unique noncommutative cone
   $C_\leq (\varphi)^+$.
   A notion of $*$-curvature $\text{ curv}_*(\varphi, D)$ 
   and Euler characteristic $\chi(\varphi, D)$ are associated with 
   each ordered pair $(\varphi, D)$, where $D\in C_\leq (\varphi)^+$.
  In this section, 
  we  obtain  asymptotic 
   formulas and basic properties for both the $*$-curvature 
    and the Euler characteristic associated with $(\varphi, D)$.  
 In the particular case when $\cH$ is a
   contractive Hilbert modules over $\CC\FF_n^+$ and $D:=I$, our two variable
    invariant
   $$
   F(\varphi, I):=(\|\varphi^*(I)\|, \text{ curv}_*(\varphi, I))
   $$ 
   is a refinement
   of the curvature invariant  from \cite{Po-curvature}  and \cite{Kr}.

Let $\CC \FF_n^+$ be the complex free semigroup algebra  generated by the free semigroup 
$\FF_n^+$ with generators $g_1,\dots, g_n$ and neutral element $e$.
Any $n$-tuple $T_1,\dots, T_n $ of bounded operators on a  Hilbert space 
$\cH$ gives rise to  a Hilbert  (left) module over  $\CC \FF_n^+$ 
in the natural way
$$
f\cdot h:= f(T_1,\dots, T_n)h, \quad f\in  \CC \FF_n^+, h\in \cH.
$$
 We associate with the  canonical operators $T_1,\dots, T_n $ 
 the completely positive linear map
 $$
 \varphi(X):= \sum_{i=1}^n T_i X T_i^*,\quad X\in B(\cH).
 $$
 The {\it adjoint} of $\varphi$ is  defined by
 $\varphi^*(X):= \sum_{i=1}^n T^*_i X T_i$. We associate 
 with  $\varphi$ and each positive  operator $D\in B(\cH)$ such that 
  $\varphi(D)\leq D$
 a  two variable numerical invariant
 \begin{equation}
 \label{F}
 F(\varphi, D):= \left(\|\varphi^*(I)\|, \,\text{\rm curv}_*(\varphi, D)\right),
 \end{equation}
 where the $*$-curvature is defined by 
 \begin{equation}\label{curv*}
 \text{\rm curv}_*(\varphi, D):=
 \lim_{k\to\infty}
  {\frac {\text{\rm trace}\,[K_{\varphi, D}^* (P_{\leq k}\otimes I)
   K_{\varphi,D}]} { 1+\|\varphi^*(I)\|+\cdots + \|\varphi^*(I)\|^k}}.
 \end{equation}
   Here, $K_{\varphi, D}$ is the Poisson 
 kernel associated with $\varphi$ and $D$ (see Section \ref{Poisson}).
 Notice that  the operator
 $K_{\varphi, D}^* (P_{\leq k}\otimes I)
   K_{\varphi,D}$ is the Poisson transform of the orthogonal  projection
   $P_{\leq k}:= I-\sum\limits_{|\alpha|=k+1} S_\alpha S_\alpha^*$.

In what follows we show that the limit defining the $*$-curvature exists.
  
 \begin{theorem}\label{curva}
 Let $\varphi$ be a w*-continuous completely positive linear map on $B(\cH)$ such that
 $\varphi(X)=\sum_{i=1}^n T_i X T_i^*$ $(n\in \NN)$. 
 If $D\in B(\cH)$ 
 is a positive operator such that  $\varphi(D)\leq D$, then
  \begin{equation}\label{curv**}
 \text{\rm curv}_*(\varphi, D):=
 \lim_{k\to\infty}
  {\frac {\text{\rm trace}\,[K_{\varphi, D}^* (P_{\leq k}\otimes I)
   K_{\varphi,D}]} { 1+\|\varphi^*(I)\|+\cdots + \|\varphi^*(I)\|^k}}
 \end{equation}
exists.
Moreover, 
~$\text{\rm curv}_*(\varphi, D)<\infty$ if and only if 
 ~$\text{\rm trace}(D-\varphi(D))<\infty$. 
\end{theorem}
\begin{proof}
Due to the properties of the Poisson kernel $K_{\varphi,D}$ 
of  Section \ref{Poisson} (see relation \eqref{r=1}), we have
\begin{equation*}\begin{split}
K_{\varphi, D}^* (P_{\leq k}\otimes I)
   K_{\varphi,D}&= 
   K_{\varphi, D}^* (I-\phi_S^{k+1} (I)\otimes I)
   K_{\varphi,D}\\
   &=
   K_{\varphi, D}^* K_{\varphi,D}- K_{\varphi,D}^*(\phi_S^{k+1} (I)\otimes I)
   K_{\varphi,D}\\
   &= D-\varphi^\infty(D)-\varphi^{k+1}(D) +\varphi^\infty (D)\\
   &= D-\varphi^{k+1}(D),
\end{split}
\end{equation*}
where $\phi_S(Y):= \sum_{i=1}^n S_iYS_i^*$ and $S_1,\ldots, S_n$ are 
the left creation operators
on the full Fock space $F^2(H_n)$.
Since the sequence $\{D-\varphi^k(D)\}_{k=1}^\infty$ is increasing,
 it is clear that 
 ~$\text{\rm curv}_*(\varphi, D)=\infty$  whenever
 ~$\text{\rm trace}(D-\varphi(D))=\infty$. Assume now that 
 ~$\text{\rm trace}(D-\varphi(D))<\infty$.
 First we consider the case when $\|\varphi^*(I)\|>1$.
Using the definition  \eqref{curv*}, we have
\begin{equation*} 
 \text{\rm curv}_*(\varphi, D)=(\|\varphi^*(I)\|-1)
 \lim_{k\to\infty}
  {\frac {\text{\rm trace}\,[ D-\varphi^k(D)]} { \|\varphi^*(I)\|^k}}.
 \end{equation*}
Let us show that this limit exists.
If $X\geq 0$ is a trace class operator, then
\begin{equation}\label{trace-ine}
\begin{split}
\text{\rm trace}\,\varphi(X)&= \text{\rm trace}\,
\sum_{i=1}^n X^{1/2} A_i^* A_i X^{1/2}=\text{\rm trace}\, 
[X^{1/2} \varphi^*(I) X^{1/2}]\\
&\leq \|\varphi^*(I)\| \text{\rm trace}\, X.
\end{split}
\end{equation}
  Since
\begin{equation}\label{dphi}
D-\varphi^{k+1}(D)= D-\varphi(D)+\varphi(D-\varphi^k(D)),\ k=1,2,\ldots,
\end{equation}
we infer that $D-\varphi^{k+1}(D)$  is a trace class operator.
 From relations \eqref{trace-ine}
and \eqref{dphi}, we obtain
\begin{equation}\label{tr-ine}
\text{\rm trace}\,[D-\varphi^{k+1}(D)]\leq \|\varphi^*(I)\|
   \text{\rm trace}\, [D-\varphi^{k}(D)]+  \text{\rm trace}\,[D-\varphi(D)].
\end{equation}
Therefore, setting
$$
a_k:= {\frac {\text{\rm trace}\, [D-\varphi^{k}(D)]} {\|\varphi^*(I)\|^k}}-
{\frac {\text{\rm trace}\, [D-\varphi^{k-1}(D)]} {\|\varphi^*(I)\|^{k-1}}},
$$
we have 
$$
a_k\leq {\frac {\text{\rm trace}\, [D-\varphi(D)]} 
{\|\varphi^*(I)\|^k}}, \quad k=1,2,\ldots.
$$
Notice that $\sum\limits_{k:\, a_k\geq 0} a_k <\infty$. Furthermore,
 every partial sum of the negative $a_k$'s is greater
 then or equal to ~$-\text{\rm trace}\, [D-\varphi(D)]$, so 
 ~$\sum_{k:\, a_k<0} a_k$ converges as well.
Therefore, 
$$
\lim_{k\to \infty}
{\frac {\text{\rm trace}\, [D-\varphi^k(D)]} 
{\|\varphi^*(I)\|^k}}
$$
exists. 
 
 Now,  assume that $\|\varphi^*(I)\|\leq 1$. Using relation \eqref{trace-ine}, we get
 $$
 0\leq  \text{\rm trace}\, [\varphi^{k+1}(D-\varphi(D))]\leq
 \|\varphi^*(I)\|  \text{\rm trace}\, [\varphi^{k}(D-\varphi(D))]
 $$
 for any $k=0,1, \ldots$.
 Hence, the sequence
  $\{\text{\rm trace}\, [\varphi^{k}(D-\varphi(D))]\}_{k=0}^\infty$
  is decreasing and
   $$
   \lim_{k\to\infty} \text{\rm trace}\, [\varphi^{k}(D-\varphi(D))]
   $$
    exists.
    If $\|\varphi^*(I)\|=1$, then, using 
      an elementary classical result and previous computations,  we obtain
    \begin{equation*}\begin{split}
 \text{\rm curv}_*(\varphi, D):&=
 \lim_{k\to\infty}
  \frac {\text{\rm trace}\,[K_{\varphi, D}^* (P_{\leq k}\otimes I)
   K_{\varphi,D}]} { n}\\
  &=\lim_{k\to\infty}
   \text{\rm trace}\,[K_{\varphi, D}^* (P_{ k}\otimes I)
   K_{\varphi,D}] \\
  &= \lim_{k\to \infty} \text{\rm trace}\, [\varphi^{k}(D-\varphi(D))],
  \end{split}
      \end{equation*}
      where $P_k:= P_{\leq k}-P_{\leq k-1}$.
 Now,  assume $\|\varphi^*(I)\|<1$.  According to the definition, we get
 $$
  \text{\rm curv}_*(\varphi, D)=(1-\|\varphi^*(I)\|) 
  \lim_{k\to \infty} \text{\rm trace}\, [D-\varphi^{k}(D)].
  $$
  Iterating relation \eqref{tr-ine}, we deduce that the sequence
  $\{\text{\rm trace}\, [D-\varphi^{k}(D)]\}_{k=1}^\infty$ is bounded.
  Since $\{D-\varphi^{k}(D)\}_{k=1}^\infty$ is increasing, we infer that
 the above  limit   exists. 
 The proof is complete.
\end{proof}

\begin{corollary}\label{curvac}
If $\varphi$ and $D$ are as in Theorem \ref{curva}, then
\begin{equation}\label{cur**}
 \text{\rm curv}_*(\varphi, D)=
 \begin{cases}
 (\|\varphi^*(I)\|-1)
 \lim\limits_{k\to\infty}
  {\frac {\text{\rm trace}\,[ D-\varphi^k(D)]} { \|\varphi^*(I)\|^k}}&
  \quad \text{ if }\  \|\varphi^*(I)\|>1,\\
  \lim\limits_{k\to \infty} \text{\rm trace}\, 
  [\varphi^{k}(D-\varphi(D))]& \quad 
  \text{ if }\  \|\varphi^*(I)\|=1,\\
  (1-\|\varphi^*(I)\|) 
  \lim\limits_{k\to \infty} \text{\rm trace}\, [D-\varphi^{k}(D)]& 
  \quad 
  \text{ if } \ \|\varphi^*(I)\|<1.
  \end{cases}
 \end{equation}
\end{corollary}

 Let $\Lambda$ be a nonempty  set of positive numbers $\alpha>0$ such that
 $$
 \text{\rm trace}\,\varphi(X)\leq \alpha \,\text{\rm trace}\,X
 $$
 for any positive trace class operator $X\in B(\cH)$. Notice that
    Theorem \ref{curva} and  Corollary
  \ref{curvac} remain true  (with exactly the same proofs)
   if we replace $\|\varphi^*(I)\|$ with  $\alpha\in \Lambda$.
   The corresponding curvature is denoted by
   $\text{\rm curv}_\alpha(\varphi, D)$.
   
   When $d:=\inf \Lambda$,  the curvature 
   $\text{\rm curv}_d(\varphi, D)$ is called the distinguished
    curvature associated with $(\varphi, D)$ and with respect to $\Lambda$.
  Now, using the analogues of Theorem \ref{curva} and  Corollary
  \ref{curvac}   for the curvatures  $\text{\rm curv}_\alpha(\varphi, D)$,
  $\alpha\in \Lambda$, one can easily prove the following.
  If  $\text{\rm curv}_\alpha(\varphi, D)>0 $, then 
  $$
   \text{\rm curv}_\alpha(\varphi, D)=
   \begin{cases}
    \text{\rm curv}_d(\varphi, D)&\quad  \text{ if } ~\alpha\geq 1,\\
    \frac{1-\alpha} {1-d}
     \,\text{\rm curv}_d(\varphi, D)&\quad  \text{ if } ~\alpha< 1.
   \end{cases}
 $$
 As we will see later in the paper,
 if $\text{\rm curv}_{\alpha_0}(\varphi, D)=0 $ for some $\alpha_0\in \Lambda$, then,
 in general, $\text{\rm curv}_d(\varphi, D)\neq 0 $. Therefore, the 
 distinguished
    curvature $\text{\rm curv}_d(\varphi, D)$ is a refinement of all the other 
    curvatures $\text{\rm curv}_\alpha(\varphi, D)$, ~$\alpha\in \Lambda$.

Let us consider an important particular case. 
According to the inequality \eqref{trace-ine}, we can take 
$\Lambda$ to be the set of all positive
 numbers $\alpha= \|\sum\limits_{i=1}^m T_i^* T_i\|$, where $(T_1, \ldots, T_m)$
  is any $m$-tuple representing the completely positive map $\varphi$, i.e.,
  $$
  \varphi(X)= \sum\limits_{i=1}^m T_i X T_i^*, \quad  X\in B(\cH).
  $$ 
  Therefore, we can talk about a distinguished
    curvature associated with $(\varphi, D)$.
All the results of this section concerning the $*$-curvature  have analogues
(and similar proofs) for the distinguished
    curvature.

From now on, for the sake of simplicity, we assume that
 $D-\varphi(D)$ is a finite
rank operator.
Using  formula \eqref{cur**}, we can  prove
 some properties of the $*$-curvature. 
 We consider  only the case when $\|\varphi^*(I)\|>1$. 
  The other cases can be treated 
  similarly, but we leave this task to the reader.

\begin{theorem}\label{curvelor}  Let $\varphi$ and $\psi$ be  w*-continuous 
completely positive linear maps on $B(\cH)$ and $B(\cH')$, respectively.
\begin{enumerate}
\item[(i)] If $X\in B(\cH)$, $Y\in B(\cH')$ are 
positive operators such that $\varphi(X)\leq X$ and $\psi(Y)\leq Y$, then
$$
\text{\rm curv}_*(\varphi\oplus \psi, X\oplus Y)=
\begin{cases}
\text{\rm curv}_*(\varphi, X)+
\text{\rm curv}_*(\psi, Y)&\quad
\text{ if } \ \|\varphi^*(I)\|=\|\psi^*(I)\|,\\ 

 \text{\rm curv}_*(\varphi, X)&\quad
\text{ if } \ \|\varphi^*(I)\|>\|\psi^*(I)\|,\\
\text{\rm curv}_*(\psi, Y)&\quad
\text{ if } \ \|\varphi^*(I)\|<\|\psi^*(I)\|.
\end{cases}
$$
\item[(ii)]
If $D_j\in B(\cH)$, $j=1,2$, are positive operators such that 
$\varphi(D_j)\leq D_j$, then
$$
\text{\rm curv}_*(\varphi, c_1 D_1+ c_2 D_2)=
 c_1\text{\rm curv}_*(\varphi, D_1)+ c_2 \text{\rm curv}_*(\varphi, D_2),
 $$
 for any positive constants $c_1,c_2$.
\item[(iii)] If $ D$  is a positive operator such that $\varphi(D)\leq D$, 
then 
$$
\text{\rm curv}_*(\varphi,D)\leq 
\text{\rm trace}\, [D-\varphi(D)]
\leq
 \|D-\varphi(D)\|\,
\text{\rm rank}\,[D-\varphi(D)]. 
$$
\end{enumerate}
\end{theorem}
\begin{proof}
According to Corollary \ref{curvac} and taking into account that
$$
\|(\varphi\oplus \psi)^*(I)\|= \max \{\|\varphi^*(I)\|, \|\psi^*(I)\|\},
$$
we have
\begin{equation*}\begin{split}
\text{\rm curv}_*(\varphi\oplus \psi, X\oplus Y)
&=
\text{\rm curv}_*(\varphi, X)\lim_{k\to \infty} \left(
\frac {\|\varphi^*(I)\|} { \|(\varphi\oplus \psi)^*(I)\|}
\right)^k \cdot\frac {\|(\varphi\oplus \psi)^*(I)\|-1}{\|\varphi^*(I)\|-1}\\
&+
\text{\rm curv}_*(\psi, Y)\lim_{k\to \infty} \left(
\frac {\|\psi^*(I)\|} { \|(\varphi\oplus \psi)^*(I)\|}
\right)^k \cdot \frac {\|(\varphi\oplus \psi)^*(I)\|-1}{\|\psi^*(I)\|-1}.
\end{split}
\end{equation*}
Hence, (i) follows. To prove (ii), notice that if $c_j\geq 0$  and 
$\varphi(D_j)\leq D_j$, then 
$$\varphi(c_1 D_1+ c_2 D_2)\leq c_1 D_1+ c_2 D_2.
$$
Taking into account Corollary \ref{curvac} and the linearity of the trace, 
 we complete the proof of (ii).
Now, using  the inequality \eqref{tr-ine}, we deduce
\begin{equation*}\begin{split}
{\frac {\text{\rm trace}\, [D-\varphi^{k+1}(D)]} 
{\|\varphi^*(I)\|^{k+1}}}
&\leq \sum_{m=1}^{k+1} {\frac {\text{\rm trace}\, [D-\varphi(D)]} 
{\|\varphi^*(I)\|^m}}\\
&=
{\frac {\text{\rm trace}\, [D-\varphi(D)]} 
{\|\varphi^*(I)\|^{k+1}}}\cdot
{\frac {\|\varphi^*(I)\|^{k+1} -1} {\|\varphi^*(I)\| -1}}.
\end{split}
\end{equation*}
Hence, and using relation \eqref{curv**}, we have
\begin{equation*}\begin{split}
\text{\rm curv}_*(\varphi,D)&\leq \text{\rm trace}\, [D-\varphi(D)]\\
&\leq
 \|D-\varphi(D)\|\,
\text{\rm rank}\,[D-\varphi(D)],
\end{split}
\end{equation*}
 and (iii) follows. The proof is complete.
\end{proof}
\begin{corollary}
 Let $\varphi$   be  a w*-continuous 
completely positive linear map on $B(\cH)$, and 
 let $D\in B(\cH)$, $D\geq 0$, be such that $\varphi(D)\leq D$.
If $D=R+Q$ is the canonical decomposition of $D$ with respect
to $\varphi$, then
$$
\text{\rm curv}_*(\varphi,D)=\text{\rm curv}_*(\varphi,Q),
$$
where $Q$ is the pure part of $D$.
\end{corollary}

Using  a result from \cite{Po-curvature}, we can prove 
the following characterization  of Hilbert modules
 isomorphic to 
  finite rank free Hilbert modules.

\begin{proposition}\label{free}
A  pure  Hilbert module $\cH$ over $\CC \FF_n^+$ is  isomorphic to 
a  finite rank free Hilbert module
$F^2(H_n)\otimes \cK$, where $\cK$ is a  Hilbert space,  if and only if  
$\varphi(I)\leq I$ and
\begin{equation}\label{fnr}
F(\varphi, I)=(n, \,\text{\rm rank}\,\cH),
\end{equation}
where $\varphi$ is the completely positive linear map associated with
the Hilbert module $\cH$.
 Moreover, in this case, $\dim \cK= \text{\rm rank}\,\cH$.
\end{proposition}

\begin{proof}
Let $T_1,\ldots, T_n$  be the canonical operators  associated with $\cH$.
 Assume that $\cH$  is isomorphic to 
a  finite rank free Hilbert module
$F^2(H_n)\otimes \cK$, i.e.,  $\cK$ is a finite dimensional Hilbert space and 
 there is a unitary operator $U:\cH\to F^2(H_n)\otimes \cK$
such 
that 
$$
T_i=U^* (S_i\otimes I_\cK) U, \quad i=1,\ldots, n.
$$
A simple calculation shows that:
\begin{enumerate}
\item[(i)]
$\text{\rm rank}\,\cH:= \dim \overline{(I-\sum_{i=1}^n T_iT_i^*)\cH}=\dim \cK$;
\item[(ii)] $\|\varphi^*(I)\|=n$;
\item[(iii)] 
$\text{\rm curv}_*(\varphi,I)=\dim \cK$.
 \end{enumerate}
 The latter equality is a consequence of Corollary \ref{curvac}.
Therefore,  relation \eqref{fnr} is satisfied.

Conversely, assume that $\cH$ is a pure  Hilbert module over $\CC \FF_n^+$
such that $\varphi(I)\leq I$ and relation  \eqref{fnr} holds.
 Then $\|\varphi^*(I)\|=n$ and 
 $\text{\rm curv}_*(\varphi,I)= \text{\rm rank}\,\cH.
 $
 Hence, we  obtain 
 $\text{\rm curv}(\cH)= \text{\rm rank}\,\cH$, where 
\begin{equation}\label{curvo}
\text{\rm curv}(\cH)=
\begin{cases} (n-1)
 \lim\limits_{k\to\infty}
  {\frac {\text{\rm trace}\,[ I-\varphi^k(I)]} 
  { n^k}}&\quad \text{ if }\  n\geq 2,\\
  \lim\limits_{k\to\infty} \text{\rm trace}\,[\varphi^k(I-\varphi(I))]
  &\quad \text{ if }\  n=1,
  \end{cases}
\end{equation} 
is the curvature invariant    introduced in  \cite{Po-curvature}.
 Now, using   Theorem 3.4 from \cite{Po-curvature}, we infer  that $\cH$ 
 is  isomorphic to 
a  finite rank free Hilbert module
$F^2(H_n)\otimes \cK$, where $\cK$ is a Hilbert space  with  
 $\dim \cK= \text{\rm rank}\,\cH$. The proof is complete.
\end{proof}

Let $\cH$ be a contractive Hilbert module over $\CC\FF_n^+$ and 
let $\varphi$ be the completely positive linear map associated with $\cH$.
What is the connection between the invariant $F(\varphi, I)$
 and the curvature invariant 
$\text{\rm curv}(\cH)$$?$
It is clear that if $\|\varphi^*(I)\|=n$,
 then  $F(\varphi, I)=(n, \,\text{\rm curv}(\cH))$.
 On the other hand,  we can prove the following result.
 
 \begin{lemma}\label{curv-curv}
 Let $\cH$ be a  finite rank contractive Hilbert module over 
 $\CC\FF_n^+$   and 
let $\varphi$ be the completely positive  linear map associated with $\cH$.
If 
  $0<\|\varphi^*(I)\|<n$, then   
  $\text{\rm curv}(\cH)=0$.
\end{lemma}
\begin{proof} First, consider the case $n\geq 2$.
Assume $\text{\rm curv}(\cH)>0$ and $\|\varphi^*(I)\|>1$.
Since $\cH$ be a  finite rank contractive Hilbert module,
 using Theorem \ref{curvelor}, we infer that
 $$
  ~\text{\rm curv}_*(\varphi,I)\leq 
\text{\rm rank}\,[I-\varphi(I)]<\infty.
$$
Taking into account relations \eqref{cur**} and \eqref{curvo}, we have
$$
{\frac { \|\varphi^*(I)\|-1)} {n-1}}
\lim_{k\to \infty}
  {\frac {\text{\rm trace}\,[ I-\varphi^k(I)]} { \|\varphi^*(I)\|^k}}
  \cdot\left[{\frac {\text{\rm trace}\,[ I-\varphi^k(I)]} { n^k}}\right]^{-1}= 
 {\frac {\text{\rm curv}_*(\varphi,I)} {\text{\rm curv}(\cH)}}<\infty.
 $$
Hence, we infer that $\lim\limits_{k\to\infty} 
{\frac{n^k} {\|\varphi^*(I)\|^k}}<\infty$, which a contradiction.
Hence, $\text{\rm curv}(\cH)=0$.
Similarly, one can show that the same conclusion holds if
 $0<\|\varphi^*(I)\|\leq 1$.
 
 Now, consider the case $n=1$. According to \cite{Po-curvature}, we have 
 $$
 \text{\rm curv}(\cH)=
 \lim\limits_{k\to\infty}
  {\frac {\text{\rm trace}\,[ I-\varphi^k(I)]} 
  { k}}
 $$
 As above, taking  again into account relation \eqref{cur**} (the case
 $\|\varphi^*(I)\|<1$) , one can  easily show that $\text{\rm curv}(\cH)=0$.
\end{proof}

Therefore, the curvature invariant  ${\text{\rm curv}(\cH)}$ does not
 distinguish  among the Hilbert modules 
 over $\CC\FF_n^+$   with  $0<\|\varphi^*(I)\|<n$.
   However, in this case,
   our $*$-curvature $\text{\rm curv}_*(\varphi, I)$ in not zero in general. 
   More precisely, we can prove the following.
   
   \begin{proposition}\label{surj}
   If $m=2,3,\ldots, n-1$, and
  $t\in (0,1]$,  then there
 exists a  finite rank contractive Hilbert module $\cH$ over $\CC \FF_n^+$, such
 that
 $\|\varphi^*(I)\|=m$ and 
  $\text{\rm curv}_*(\varphi, I)=t$, i.e.,
 $$
 F(\varphi, I)= (m,t).
 $$
   \end{proposition}
 \begin{proof}
 According to Theorem 3.8 from \cite{Po-curvature}, there is 
  a  finite rank contractive Hilbert module $\cH$ over $\CC \FF_m^+$
  such that $\text{\rm curv}(\cH)=t$.
  Let $T:=[T_1,\ldots, T_m]$ be the row contraction associated with $\cH$ 
  and let $\varphi_T$ be the corresponding
  completely positive map. Since $ t>0$, Lemma \ref{curv-curv} implies 
  $\|\varphi^*(I)\|=m$.
  Let $\cH'$ be the 
  finite rank contractive Hilbert module over $\CC \FF_n^+$ defined by the row 
  contraction
  $$
  [T_1, \ldots, T_{m-1}, 
  \frac {1} {\sqrt{n-m+1}} T_m, \ldots, \frac {1} {\sqrt{n-m+1}} T_m],
 $$
 and let $\varphi$ be the associated completely positive map.
 Notice that $\varphi(X)= \varphi_T(X)$, $X\in B(\cH)$, and 
 $\varphi^*(I)= \varphi_T^*(I)$. Hence, we have $\|\varphi^*(I)\|=m$.
 Taking into account  relations \eqref{cur**} and \eqref{curvo}, we infer that
 \begin{equation*}
 \begin{split}
 \text{\rm curv}_*(\varphi, I)&= 
 (m-1)
 \lim_{k\to\infty}
  {\frac {\text{\rm trace}\,[ I-\varphi_T^k(I)]} { m^k}}\\
  &= \text{\rm curv}(\cH)=t.
 \end{split}
 \end{equation*}
 Summing up, we have $F(\varphi, I)=(m,t)$, which completes the proof.
 \end{proof}
 This result clearly  shows that, in the particular case of  finite rank 
 contractive 
Hilbert modules over $\CC\FF_n^+$,
our   invariant  $F(\varphi, I)$ is a refinement of  $\text{\rm curv}(\cH)$.

\bigskip

Let $\cH$ be a Hilbert module over $\CC\FF_n^+$, ~ $\varphi$ be its 
associated completely positive map,
and  let $ D $ be a positive operator   such that $\varphi(D)\leq D$.
For each $k=0,1,\ldots,$ define
$$
M_k(\varphi, D):=\text{\rm span} \left\{p\cdot \xi:\ p\in \CC\FF_n^+, 
~\text{\rm deg}(p)\leq k, \ \xi\in [D-\varphi(D)]^{1/2} \cH\right\}.
$$ 
   We define the Euler characteristic associated with $\varphi$ and $ D$ 
   by setting
\begin{equation}\label{chi}
\chi(\varphi, D):= \lim_{k\to\infty} 
{\frac { \text{\rm dim} \,M_k(\varphi,D)}{1+n+\cdots + n^k}}.
\end{equation}
In what follows we show that the limit \eqref{chi} exists (finite or infinite).

\begin{theorem}\label{euler}
Let $\varphi$ be a w*-continuous completely positive linear map on $B(\cH)$ such that
 $\varphi(X)=\sum_{i=1}^n T_i X T_i^*$ $(n \geq 2)$, and 
let $ D$  be a positive operator    such that $\varphi(D)\leq D$. 
Then the Euler characteristic $\chi(\varphi, D)$ exists and 
\begin{equation*}
\begin{split}
\chi(\varphi, D)
&=
\lim_{k\to\infty}
  {\frac {\text{\rm rank}\,[K_{\varphi, D}^* (P_{\leq k}\otimes I)
   K_{\varphi,D}]} { 1+n+\cdots + n^k}}\\
   &= (n-1)
   \lim_{k\to\infty}
  {\frac {\text{\rm rank}\,[D-\varphi^k(D)]} {  n^k}}.
\end{split}
\end{equation*}
Moreover, the Euler characteristic $\chi(\varphi, D)<\infty$ if
and only if  $\text{\rm rank}\,[D-\varphi(D)]<\infty$. 
\end{theorem}
\begin{proof}
 Notice   that if 
 $\text{\rm rank}\,[D-\varphi(D)]=\infty$,
then the inclusion $M_0(\varphi,D)\subset M_k(\varphi,D)$ implies 
 $\chi(\varphi, D)=\infty$.
Now, assume  that $\text{\rm rank}\,[D-\varphi(D)]<\infty$.
Notice that $M_k(\varphi,D)$ is equal to the range of
 $K_{\varphi, D}^*(P_{\leq k}\otimes I)$. Since the
  latter operator has finite rank, we have
  $$
  \dim M_k(\varphi,D)= \text{\rm rank}\,[K_{\varphi, D}^* 
  (P_{\leq k}\otimes I)
   K_{\varphi,D}]=\text{\rm rank}\,[D-\varphi^{k+1}(D)].
   $$  
For each $k\geq 1$, we have
$$
M_k(\varphi,D)=
M_0(\varphi,D)+ T_1M_{k-1}(\varphi,D)+ \cdots +T_nM_{k-1}(\varphi,D).
$$
Hence, we infer
$$
\dim M_k(\varphi,D)\leq n \dim M_{k-1}(\varphi,D) + \text{\rm rank}\, 
[D-\varphi(D)], \quad k=1,2,\ldots.
$$
The rest of the proof is similar to that 
of Theorem 4.1 from \cite{Po-curvature}.
We shall  omit it.
\end{proof}
We should mention that 
in the particular case when $\varphi (I)\leq I$ and $D=I$,
the result of 
 Theorem \ref{euler}
was obtain  in \cite{Po-curvature}.

Similarly to the proof of Theorem \ref{curvelor}, one can prove 
the following result.
\begin{proposition}
 For each $j=1,2$,
let $\cH_j$ be a Hilbert module over $\CC\FF_n^+$,  and  $\varphi_j$ 
be its 
associated completely positive linear map.
Let $ D_j$ be  positive operators  such that $\varphi_j(D)\leq D_j$.
Then we have:
\begin{enumerate}
\item[(i)]
$\chi(\varphi_1\oplus \varphi_2, D_1\oplus D_2)= \chi(\varphi_1, D_1)+
\chi(\varphi_2, D_2).
$
\item[(ii)]
If $ D_1\leq I$ and $\|\varphi_{1}^*(I)\|=n$, then 
$\text{\rm curv}_*(\varphi_1,D_1)\leq \chi(\varphi_1, D_1)$.
\end{enumerate}
\end{proposition}



\end{document}